


\input amstex
\documentstyle{amsppt}

\magnification=1200
\parskip 12pt
\pagewidth{5.4in}
\pageheight{7.2in}

\baselineskip=12pt
\expandafter\redefine\csname logo\string@\endcsname{}
\NoBlackBoxes                
\NoRunningHeads
\redefine\no{\noindent}

\define\C{\bold C}
\define\R{\bold R}
\define\Z{\bold Z}

\redefine\P{\bold P}

\define\al{\alpha}

\define\de{\delta}
\define\la{\lambda}
\define\si{\sigma} 
\define\Si{\Sigma} 
\define\La{\Lambda}

\define\om{\omega}

\define\sub{\subseteq}  
\define\es{\emptyset}

\def\pr{\prime} 
\define\st{\ \vert\ }   

\redefine\ll{\lq\lq}
\redefine\rr{\rq\rq\ }
\define\rrr{\rq\rq}   

\redefine\det{\operatorname {det}}
\redefine\deg{\operatorname {deg}}
\define\Hol{\operatorname {Hol}}

\redefine\dim{\operatorname {dim}}
\define\codim{\operatorname {codim}}
\define\End{\operatorname {End}}
\define\rank{\operatorname {rank}}

\define\Ker{\operatorname {Ker}}

\define\PD{\operatorname {PD}}

\define\grkcn{Gr_k(\C^n)}
\define\glnc{GL_n\C}

\define\lan{\langle}
\define\ran{\rangle}

\redefine\i{{\ssize\sqrt{-1}}\,}

\define\sla{ {}^{\la} }
\define\na{\nabla}
\define\bd{\partial}
\define\ABC{\lan A \vert B \vert C \ran}

\define\subD{_{\ssize D}}
\define\HolD{\operatorname {Hol}\subD}

\topmatter
\title Introduction to homological geometry: part I
\endtitle
\author Martin A. Guest
\endauthor
\endtopmatter

\document
Given a space $M$, one may  attempt to construct various natural 
cohomology algebras such as the ordinary (simplicial, singular, etc.)
cohomology algebra $H^\ast(M)$ and
the  quantum cohomology algebras  $QH^\ast(M)$ and $\widetilde{QH}^\ast(M)$.  
For example, if $M$ is the $n$-dimensional complex projective
space $\C P^{n}$, then
$$
H^\ast(\C P^{n}) \cong \C[p]/\lan p^{n+1}\ran,
\quad
\widetilde{QH}^\ast(\C P^{n}) \cong \C[p,q]/\lan p^{n+1}-q\ran.
$$
In the first case the algebraic variety defined by the equation
$p^{n+1}=0$ is not very interesting,
but in the second case we have the nontrivial variety
$p^{n+1}-q=0$ in $\C^2$.
Roughly speaking, \ll homological geometry\rr is concerned with
\ll geometry of the algebraic
variety $V_M$, where the algebra of functions on $V_M$ 
is the quantum cohomology algebra of the space
$M$\rrr.  

Our main sources of inspiration for this subject are the
stimulating papers \cite{Gi-Ki} and \cite{Gi1}-\cite{Gi6}.  
The excellent (if idiosyncratic) survey papers 
\cite{Au1}-\cite{Au4} amplify and explain some of this material.
Our aim in these lectures is very modest: if Audin is an introduction
to Givental, then the first part of this survey
of homological geometry will be an introduction to Audin.

We shall not discuss the rigorous definition of  quantum cohomology.
For this, the reader should consult the articles
\cite{Ko}, \cite{Ru-Ti} and the books and survey articles listed later on.
Nor do we discuss the historical motivation, which comes from physics.
Nevertheless, we shall begin by giving an informal
introduction to quantum cohomology, together with some very explicit
calculations in the appendices.  Readers who have found
quantum cohomology intimidating may wonder how it will be possible 
to contemplate applications of the theory, after such a superficial 
treatment of the foundations. My answer would be that this
is  already standard practice for {\it ordinary} cohomology
theory; the rigorous foundations of cohomology theory
are unavoidably messy, yet there is no difficulty in
computing (for example) 
the cohomology rings for simple spaces such as surfaces
or projective spaces. In quantum cohomology we face the same situation.
We define a product operation by intersecting certain special
kinds of cycles, and, even though  
the general definition is complicated, we can perform the calculations
satisfactorily for nice spaces (such as homogeneous K\"ahler manifolds). 

Our main purpose is to describe a path from this naive intersection-theoretic
formulation of quantum cohomology to  the way in which differential
equations --- especially those related to differential geometry
and the theory of integrable systems --- enter into quantum cohomology. 
We begin, therefore, with some generalities on flat connections and 
Frobenius manifolds, and then discuss quantum cohomology as an example
from this point of view.  In sections 4, 5 and 6 we give three examples of how
differential equations arise. In the last section we indicate very briefly
how the flat connection defined by quantum cohomology underlies these
phenomena.  In part II we shall consider the differential equations in more depth.

In addition to the articles mentioned earlier, the books and survey
articles \cite{BCPP}, \cite{Co-Ka}, \cite{Du}, 
\cite{Ma}, \cite{Mc-Sa}, \cite{Si}, \cite{Ti},
\cite{Vo} may be consulted for further information and many more references. 
Although we shall not reach the point of discussing any of the famous enumerative
results (because we use only the \ll small\rr quantum
product) or mirror symmetry, we hope that the reader will be prepared to go on to
these topics after reading this introduction.

{\eightpoint
Acknowledgements:  This article is based on lectures given at the Center
for Theoretical Sciences, National Tsing Hua University, Taiwan,
in March 1999.  The author is grateful to Professor Chuu-Lian Terng for her invitation
to give those lectures, and to the NCTS for its hospitality.
He also thanks Augustin-Liviu Mare, Takashi Otofuji, and Ken-ichi Seimiya
for several suggestions.}

\head
\S 1 Flat connections and Frobenius algebras
\endhead

References:  \cite{Du}, \cite{Mc-Sa}

Let $W$ be a finite-dimensional vector space (real in this section,
but complex in future sections).  Let $d$ denote
the directional derivative operator on functions $W\to W$; thus,
for $X\in W$ and $Y:W\to W$, we denote by $d_X Y$ the
derivative of $Y$ in the direction of $X$.  
We regard $d$ as the standard flat connection
(covariant derivative operator)  on the manifold $W$. 

Any other connection on $W$ is of the form $\nabla = d + \om$, where
$\om$ is a $1$-form on $W$ with values in the Lie algebra $\End(W)$.  
Given any
such $\om$, we define an associated family of connections $\na\sla$ by
$\na\sla = d + \la \om$, where $\la\in \R$.

It is easy to check that
$$
\align
\na \ \text{is flat}&\iff d\om+\om\wedge \om =0\tag 1\\
\na\sla \ \text{is flat for all $\la$}&\iff d\om=0, \om\wedge \om =0
\endalign
$$
and
$$
\align
\na \ \text{has zero torsion}&\iff \om(X)Y = 
\om(Y)X \ \text{for all}\ X,Y \tag 2\\
&\iff\na\sla \ \text{has zero torsion for all}\ \la.
\endalign
$$
Similarly, if $(\ ,\ )$ is a symmetric bilinear form on $W$, then
$$
\align
\na \ \text{is compatible with $(\ ,\ )$}&\iff
( \om(X)Y,Z) + ( Y,\om(X)Z) = 0\ \text{for all}\ X,Y,Z \tag 3\\
&\iff\na\sla 
\ \text{is compatible with $(\ ,\ )$ for all}\ \la.
\endalign
$$

We shall be interested in the special case where $W$ is a 
{\it Frobenius algebra.}
This means that $W$ has a commutative 
associative product operation $\circ$, so that
$W$ becomes an algebra, and that the bilinear form is nondegenerate
and also satisfies the \ll Frobenius condition\rr
$( X\circ Y, Z) = ( X,Y\circ Z)$ for all $X,Y,Z\in W$.

For a Frobenius algebra $W$, we shall consider the
$1$-form  $\om(X)Y=X\circ Y$, and the corresponding connections
$$
\na\sla_X Y = d_XY + \la X\circ Y.
$$

\proclaim{Proposition} For any Frobenius algebra we have
$$\align
& \na\sla \ \text{is flat for all}\ \la \tag i\\
 & \na\sla \ \text{has zero torsion for all}\ \la \tag ii\\
\endalign
$$
\endproclaim

\demo{Proof} We use conditions (1) and (2) above.  Commutativity of
$\circ$ gives (ii).  For (i) , the fact that $\om_t(X)(Y)$ is independent of $t\in W$
gives $d\om=0$, and commutativity and associativity
of $\circ$ gives $\om\wedge\om=0$.
\enddemo

The connection $\na\sla$ is compatible with the bilinear form $(\ ,\ )$
if and only if $\la=0$, i.e. $\na\sla=d$.  Indeed, the Frobenius condition
may be interpreted as a kind of \ll skew-compatibility\rrr.

More generally, one can consider a family $\circ_t$ 
of Frobenius algebra structures on $W$.  We shall be interested in
the case where the family is parametrized by $t\in W$.  Since the tangent bundle of
the manifold $W$ is trivial, we may then regard $\circ_t$  as
a Frobenius algebra structure on the tangent space
$T_tW$ at $t$. The formula
$$
\na\sla_X Y = d_XY + \la X\circ_t Y
$$
defines a connection on the manifold $W$; it is called
the {\it Dubrovin connection.} 

From the definition of Frobenius structure
it can be shown (as in the proposition) that $\na\sla$ has zero torsion,
and also that $\om\wedge\om=0$,  where $\om_t(X)(Y) = X\circ_t Y$.
However, it is not
necessarily flat, because we do not necessarily have $d\om=0$.

We are now in a situation  similar to that of
having an almost complex structure on a manifold,
i.e. where there is a complex structure on each tangent space, but these complex
structures are not necessarily \ll integrable\rrr.  The manifold $W$, with
the family of Frobenius structures $\{ \circ_t \st t\in W\}$, is in
fact an example of a pre-Frobenius manifold (see \cite{Du}).  It will be a
Frobenius manifold --- i.e. it will be regarded as \ll integrable\rr --- if the
Dubrovin connection is flat for all $\la$, i.e. if $d\om = 0$.

\head
\S 2 An example of a Frobenius algebra:  $(H^\ast(M),\circ_t)$
\endhead

Reference: \cite{Mc-Sa}

Let $M$ be a simply connected (and compact, connected) K\"ahler manifold,
of complex dimension $n$.  We shall assume that the integral cohomology
of $M$ is even dimensional and torsion-free, i.e. that

\no (A1) \quad
$H^*(M;\Z)\ =\ \bigoplus_{i=0}^n H^{2i}(M;\Z) \ \cong\ 
\bigoplus_{i=0}^n \Z^{m_{2i}}$

\no where $m_{2i}=\rank H^{2i}(M;\Z)$.

The role of the vector space of \S 1 will be played by
$W=H^*(M;\C)\cong \bigoplus_{i=0}^n \C^{m_{2i}}$.  There is a family 
 $\{\circ_t\st t\in W\}$ of Frobenius structures 
on the vector space $W$, called the
{\it quantum product.}  In this section we shall define  $\circ_t$
in terms of certain {\it Gromov-Witten invariants} $\ABC\subD$. In the
next section we shall define $\ABC\subD$ and give some examples.

We need some notation from ordinary cohomology theory.  Let
$$
\PD:H^i(M;\Z) \to H_{2n-i}(M;\Z)
$$
be the Poincar\acuteaccent e duality isomorphism.  This may be
defined as the map which sends a cohomology class 
$x$ to the \ll cap product\rr of
$x$ with the fundamental class of $M$. As far as
possible we shall use the notation
$$
a,b,c,\dots \ \in\  H^\ast(M;\Z) \quad (\  \sub H^\ast(M;\C) \  )
$$
for cohomology classes, and we write 
$\vert a\vert,\vert b\vert,\vert c\vert,\dots$
for the degrees of $a,b,c,\dots$.  We shall write
$$
A=\PD(a), B=\PD(b), C=\PD(c),\dots \ \in\  
H_\ast(M;\Z) \quad (\ \sub H_\ast(M;\C)\  )
$$
for the Poincar\acuteaccent e dual homology classes, and
$\vert A\vert,\vert B\vert,\vert C\vert,\dots$
for their degrees.

Let 
$$
\lan\ ,\ \ran:H^i(M;\Z) \times H_i(M;\Z) \to \Z
$$
denote  the Kronecker (or \ll evaluation\rrr) pairing; we use the
same notation for the extended pairing
$$
\lan\ ,\ \ran:H^\ast(M;\Z) \times H_\ast(M;\Z) \to \Z
$$
(thus, $\lan a,B \ran$ is zero whenever $\vert a\vert \ne \vert B\vert$).
Since there is no torsion, both these pairings are nondegenerate.

The role of the bilinear form  of \S 1 will
be played by the ($\C$-linear extension of the)
\ll intersection pairing\rrr, which is defined by
$$
(\ ,\ ):H^\ast(M;\Z) \times H^\ast(M;\Z) \to \Z,\quad
(a,b)=\lan ab, M\ran.
$$
On the right hand side, 
$M$ denotes the fundamental homology class of the manifold $M$.  It is
an element of $H_{2n}(M;\Z)$, and its Poincar\acuteaccent e dual
cohomology class --- the identity element of the
cohomology algebra --- will be denoted by $1\in H^{0}(M;\Z)$
(an exception to our notational convention for cohomology
classes!).  The homology class represented by a single point of
$M$ will be denoted by $Z\in H_0(M;\Z)$, and its Poincar\acuteaccent e dual
cohomology class will be denoted in accordance with our convention
by $z$.

Now, the well known duality between the cup product and
the Kronecker pairing may be expressed by the formula
$$
\lan ab, M\ran = \lan a,B\ran = \lan b,A \ran.
$$
From the nondegeneracy of the Kronecker pairing (in this
situation), it follows that the intersection pairing $(\ ,\ )$
is nondegerate.  The formula also shows that the cup product
satisfies the Frobenius condition $(ab,c)=(a,bc)$,
because both sides of this equation are equal to
$\lan abc,M\ran$.  Hence the cup product and the intersection
pairing give rise to a Frobenius algebra structure on the vector space $W$.
It should be noted that the intersection pairing is always indefinite (unless
$M$ is zero-dimensional).

The resulting Frobenius structure on the manifold $W$ is trivial,
in the sense that the associated $1$-form $\om_t(x)(y) = xy$
is constant, i.e. independent of $t$. The condition $d\om=0$ is therefore
automatically satisfied. On the other hand, the quantum product will
give a nontrivial example, as we shall see next.

The cup product on $H^\ast(M;\Z)$ may be specified in terms
of its \ll structure constants\rrr.  To do this, we choose
generators as follows:
$$
\align
H_\ast(M;\Z) &= \bigoplus_{i=0}^s \Z A_i \\
H^\ast(M;\Z) &= \bigoplus_{i=0}^s \Z a_i
\endalign
$$
and we define \ll Kronecker dual\rr cohomology classes 
$a^c_0,\dots,a^c_s$ (i.e. the dual basis
with respect to $(\ ,\ )$)  by $a_i a^c_j = \de_{ij}z$.
Then for any $i,j$ we have
$$
a_ia_j= \sum_{\{\al \st \vert a_{\al}\vert = 
\vert a_{i}\vert + \vert a_{j}\vert\} } \la_{\al}^{i,j} a^c_{\al}
$$
for some $\la_1^{i,j},\dots,\la_s^{i,j}\in \Z$.  These structure
constants are given by
$$
\la_k^{i,j}= \lan a_ia_ja_k, M\ran.
$$
Observe that if $\la^{i,j}_k\ne 0$ then the numerical condition
$$
\vert a_{i}\vert + \vert a_{j}\vert + \vert a_{k}\vert = 2n
$$
must be satisfied.  

\proclaim{Definition} For cohomology classes $a,b,c$ we define
$\ABC_0 = \lan ab,C\ran = \lan abc, M\ran$.
\endproclaim

\no Giving all the structure constants is the same thing as giving
all \ll triple products\rr $\ABC_0$.

The quantum product will be determined by a larger family of
triple products denoted by
$\ABC\subD$, where $D$ varies in $H_2(M;\Z)$.  
(These are the Gromov-Witten invariants 
$\psi^D_{3,0}(A\otimes B\otimes C)$.)
We shall define these new triple products in the next section; for the
rest of this section we shall just assume that

\no (A2)\quad
$\ABC\subD\in\Z$ is defined for any $A,B,C\in H_\ast(M;\Z)$, 
$D\in H_2(M;\Z)$.  Moreover, $\ABC\subD\in\Z$ is linear and symmetrical
in $A,B,C$.

\no To define $a\circ_t b$ for $a,b\in W$ and $t\in W$, it suffices
to define $\lan a\circ_t b, C\ran$ for all $C\in H_\ast(M;\C)$.  
We shall give a definition only for $t\in H^2(M;\C)$; the family
$\{\circ_t \st t\in H^2(M;\C) \}$ is called the \ll small quantum
product\rrr. The definition is:

\proclaim{Definition}  $\lan a\circ_t b, C\ran=
\sum_{D\in H_2(M;\Z)}\,\ABC\subD\, e^{\lan t,D\ran}$.
\endproclaim

\no We assume that

\no (A3) in the above definition, the sum on the right hand side is finite.

\no Observe that as \ll$t\to -\infty$\rr the right hand side
converges to $\ABC_0$; hence $a\circ_t b$ converges to the
cup product $ab$.  In this sense, the quantum product is a deformation
of the cup product.

The main result concerning this quantum product is:

\proclaim{Theorem}  $\{ \circ_t \st t\in W \}$ defines a
Frobenius structure on $W=H^\ast(M;\C)$.
\endproclaim

\no The proof has two nontrivial ingredients: the
associativity of the quantum product (which we shall take for
granted), and the condition $d \om=0$ (which we shall discuss in \S 7
in the case $t\in H^2(M;\C)$).

There is a modification of the quantum product, which is a product
operation
$$
\circ:H^\ast(M;\C) \times H^\ast(M;\C) \to H^\ast(M;\C)\otimes \La
$$
where $\La$ is the group algebra $\C[H_2(M;\Z)]$.  Formally, an
element of $\La$ is a finite sum $\sum_{i}\la_i q^{D_i}$,
where $\la_i\in\C$, $D_i\in H_2(M;\Z)$, and where the symbols $q^D$ are multiplied
in the obvious way, i.e. $q^Dq^E=q^{D+E}$.  (For the time
being, therefore, the symbol $q$ has special status as a
\ll formal variable\rr --- it is certainly not a cohomology class!)
The definition is:

\proclaim{Definition}  $a\circ b=
\sum_{D\in H_2(M;\Z)}\,(a\circ b)\subD\, q^D$, where $(a\circ b)\subD$
is defined by $\lan (a\circ b)\subD, C\ran=\ABC\subD$ for all
$C\in H_\ast(M;\C)$.
\endproclaim

\no We shall also refer to this operation as the quantum
product.  By our assumption (A3), the sum on the right hand side of
the definition of $a\circ b$ is finite.  

Extending the (small) quantum product in a $\La$-linear fashion,
we can make $W\otimes \La$ into an algebra.  This is called
the {\it (small) quantum cohomology algebra} and it is denoted
by $QH^\ast(M;\C)$.
The Frobenius algebra $(W,\circ_t)$ may be obtained from the
quantum cohomology algebra $QH^\ast(M;\C)$ 
by \ll putting $q^D=e^{\lan t,D\ran}$\rrr.
Thus, $(W,\circ_t)$ and $QH^\ast(M;\C)$ contain essentially
the same information.

It is convenient to define a grading on the algebra $W\otimes\La$ by
defining
$$
\vert aq^D\vert = \vert a\vert + 2\lan c_1 TM, D\ran.
$$
We shall assume that the quantum products preserve this grading,
i.e. that $\vert a\circ b\vert = \vert a\vert + \vert b\vert$.
It follows that $\vert (a\circ b)\subD\vert = \vert a\vert
+\vert b\vert - 2\lan c_1 TM,D\ran$.  Hence, if $\ABC\subD\ne 0$,
then the numerical condition $\vert a\vert + \vert b\vert
- 2\lan c_1 TM, D\ran = \vert C\vert$ must be satisfied, i.e.
$$
\vert a\vert +\vert b\vert +\vert c\vert = 2n + 2\lan c_1 TM,D\ran.
$$
The geometrical meaning of this numerical condition will 
become clear in the next section.
We write $QH^i(M;\C)=\{ x\in QH^\ast(M;\C) \st \vert x\vert = i \}$.

One further piece of notation will be useful.  In all of our examples
we shall choose an identification $H_2(M;\Z)\cong \Z^r$ (for
some $r\ge 1$).  Having made this choice, we write $D=(s_1,\dots,s_r)$,
and $q^D=q_1^{s_1}\dots q_r^{s_r}$.  In our examples --- although
not necessarily for more general manifolds $M$ --- it will turn out that the
subset $\widetilde{QH}^\ast(M;\C) = 
H^\ast(M;\C) \otimes \C[q_1,\dots,q_r]$ is actually a {\it subalgebra} of 
$QH^\ast(M;\C)$ (with respect to the quantum product).  
This implies in particular that
$QH^\ast(M;\C)$ and $\widetilde{QH}^\ast(M;\C)$ contain 
the same information, and so we can restrict attention to
$\widetilde{QH}^\ast(M;\C)$.

\head
\S 3  Examples of $(H^\ast(M),\circ_t)$ and $QH^\ast(M)$
\endhead

References: \cite{Gi-Ki}, \cite{Mc-Sa}, \cite{Fu-Pa}

Quantum cohomology depends on the triple products
$\ABC\subD$, and we cannot postpone their definition any longer.  
In this section we shall give the \ll naive\rr definition of
$\ABC\subD$ (as in \cite{Gi-Ki}, for example); this has the advantage of
simplicity, but the disadvantage that it will be impossible to
prove any general theorems.   Our main objective is simply to
describe some concrete examples as preparation for homological geometry.
In this section we shall concentrate on the example $M=\C P^n$, but
later on we consider $M=\grkcn$ (the Grassmannian of complex $k$-planes in $\C^n$
--- in \S 4), $M=F_n$ (the space of complete flags in $\C^n$ --- in
Appendix 1), and $M=\Si_k$ (the Hirzebruch surface --- in Appendix 2).

It has been said that quantum cohomology is \ll easy to
calculate, but hard to define\rrr.  This is not as strange as it sounds,
because exactly the same can be said for {\it ordinary} cohomology,
if one restricts attention to a few nice spaces.
In this spirit, we shall begin by reviewing the definition of the
triple product $\ABC_0 = \lan ab,C\ran = \lan abc, M\ran$ in ordinary cohomology.

The naive definition is
$$
\ABC_0 \ =\ \vert \tilde A\cap \tilde B \cap \tilde C \vert,
$$
where the right hand side means the number of points (counted
with multiplicity) in the intersection $\tilde A\cap \tilde B \cap \tilde C$,
where $\tilde A,\tilde B,\tilde C$ are suitable representatives
of the homology classes $A,B,C$.  In certain situations this definition
is \ll correct\rrr, in the sense that it gives the usual triple product
$\lan abc, M\ran$.   (Similarly, $\lan ab, M\ran $ may be defined naively as
$\vert \tilde A\cap \tilde B\vert$.) For example, in the complex algebraic
category, the definition is correct whenever there exist
representative algebraic subvarieties $\tilde A,\tilde B,\tilde C$
whose intersection is finite (or empty) --- see the appendix
of \cite{Fu}.  Thus, whenever we are lucky enough to find such
representatives, we can calculate $\ABC_0$.  

The most famous
example where this method works is the case $M=\grkcn$
(Schubert calculus).  Here all the 
generators of the homology groups are representable by algebraic
cycles (Schubert varieties), and for any three such generators
$a,b,c$ satisfying the condition 
$\vert a\vert + \vert b \vert + \vert c\vert = \dim M$ 
there exist representatives whose intersection is finite (or empty).

\proclaim{Example: $M=\C P^n$}
\endproclaim

We have $H^{2i}(\C P^n;\Z)\cong \Z x_i$ $(0\le i\le n)$,
where the Poincar\acuteaccent e dual homology generator $X_i$
(of degree $2n-2i$) can be represented by $\tilde X_i = \P(V)$,
for any complex linear subspace $V\sub \C^{n+1}$ of 
codimension $i$.  Following our usual notational conventions,
we shall write $x_0=1$ and $x_n=z$.

If $i+j>n$, then $x_ix_j=0$, since $H^{i+j}(\C P^n;\Z)=0$.
If $i+j\le n$, then $x_ix_j=\la x_{i+j}$ for some $\la$,
and we have to calculate $\la = \lan x_ix_jx^c_{i+j}, \C P^n \ran=
\lan X_i \vert X_j \vert X^c_{i+j} \ran_0$.

First, we have
$x_ix_{n-j}=\de_{ij}z$, as there exist linear
subspaces $V,W$ of $\C^{n+1}$ of codimensions $i,n-j$
such that $\P(V)\cap \P(W)$ is finite and nonempty
if and only if $i=j$, and in this case the intersection
consists of a single point (of multiplicity one).  This
shows that $x^c_i=x_{n-i}$.

To calculate $\lan X_i \vert X_j \vert X^c_{i+j} \ran_0
= \lan X_i \vert X_j \vert X_{n-(i+j)} \ran_0$
(for $i+j\le n$), we represent the three classes by
linear subspaces of $\C^{n+1}$ of codimensions $i,j,n-(i+j)$.
If the subspaces are in general position, the codimension of the
intersection is $i+j+n-(i+j)=n$, so this triple intersection
is a line, and $\tilde X_i \cap \tilde X_j \cap \tilde X_{n-(i+j)}$
is a single point of $\C P^n$.  Since we are taking intersections of
linear subspaces, the multiplicity of this point is one. We
conclude that $\la = 1$, and so $x_ix_j=x_{i+j}$.

The cohomology algebra of $\C P^n$ is therefore isomorphic
to $\C[p]/\lan p^{n+1}\ran$ (where $p=x_1$), as stated in the introduction.

We wish to define $\ABC\subD$ as a certain intersection number.
From our discussion in the last section, it
should have the properties

\no(a)\quad
$D=0 \implies \ABC\subD=\lan abc, M\ran$

\no(b)\quad $\ABC\subD\ne 0\implies 
\vert a\vert +\vert b\vert +\vert c\vert = 2n + 2\lan c_1 TM,D\ran$.

\no (Actually, (a) will be obvious from the definition, but (b) will have
the status of an assumption.)

\no The naive definition is
$$
\ABC\subD=\vert
\HolD^{\tilde A,p} \cap \HolD^{\tilde B,q} \cap \HolD^{\tilde C,r} 
\vert
$$
where
$$
\HolD^{\tilde A,p}=
\{ \text{holomorphic maps}\ f:\C P^1\to M \st
f(p)\in\tilde A\ \text{and}\  [f]=D \}
$$
and where $\tilde A$ is a representative of the homology
class A.  The points $p,q,r$ are three distinct basepoints in $\C P^1$.
The notation $[f]$ denotes the homotopy class of $f$, which is an
element of $\pi_2(M)\cong H_2(M;\Z)$.

We {\it assume} that there exist $\tilde A,\tilde B, \tilde C$ such
that the above intersection is finite (or empty); this is essentially
the meaning of assumption (A2) of the previous section.  More
fundamentally, we shall make the following three assumptions:

\no(A2a)\quad
$\HolD^{M,p}$ is a (smooth) complex manifold, of complex dimension
$n+\lan c_1TM,D\ran$. This will be denoted more briefly by $\HolD$.

\no(A2b)\quad
$\HolD^{\tilde A,p}$ is a  complex submanifold
(or subvariety)  of $\HolD^{M,p}$, and
the complex codimension of $\HolD^{\tilde A,p}$ in $\HolD^{M,p}$
is equal to the complex codimension of $\tilde A$ in $M$.

\no(A2c)\quad
If $\tilde A,\tilde B, \tilde C$ intersect transversely, then so
do $\HolD^{\tilde A,p}, \HolD^{\tilde B,q}, \HolD^{\tilde C,r}$.

\no We shall not discuss the extent to which (A2) is equivalent
to (A2a), (A2b), (A2c), nor the extent to which these conditions
are true. (Certainly all four assumptions hold for homogeneous
K\"ahler manifolds.)  However, a brief comment on
the origin of the numerical expressions in (A2a) and (A2b) 
may be helpful.
First, the tangent space at $f\in \HolD^{M,p}$ --- assuming that
$\HolD^{M,p}$ is a manifold --- may be identified with the
space of holomorphic sections of the bundle $f^\ast TM$.  By
the Riemann-Roch theorem, the complex dimension of this vector space is
$n+\lan c_1TM,D\ran$, as stated in (A2a).  Regarding (A2b), this
would follow from the commutative diagram
$$
\CD
\HolD^{\tilde A,p} @>\sub>> \HolD^{M,p}\\
@VVV @VVV\\
\tilde A @>\sub>> M
\endCD
$$
in which the vertical maps are given by $f\mapsto f(p)$, whenever
the evaluation map $\HolD^{M,p}\to M$ is a locally trivial fibre bundle
(as it will be when $M$ is homogeneous, for example).

It is important to bear in mind that the space $\HolD$ is
finite-dimensional and algebraic, so there is no question of
infinite-dimensional analysis here.  The technical problems in
giving a rigorous definition arise from the noncompactness of $\HolD$,
and the fact that the transversality condition may not be satisfied.
At the end of this section we shall comment briefly on the kind of technical
conditions that are needed.

To understand the definition better, let us consider the case
where $\pi_2(M) \cong \Z$, e.g. $M=\C P^n$ or $\grkcn$.  Then we
may write $D = sD_0$ where $s\in \Z$ and $D_0$ is a generator
of $H_2(M;\Z)$.  For a map $f:\C P^1\to M$,
we have $[f] = sD_0$ for some $s$, and we consider this $s$
to be the \ll degree\rr of $f$.  

Let $\lan c_1 TM,D_0\ran = N$.  We shall assume here that $N>0$
(we have $N=n+1$ for $\C P^n$ and $N=n$ for $\grkcn$).
Then if $\ABC_s$, i.e. $\ABC\subD$,  is nonzero, we have the numerical
condition
$$
\vert a\vert +\vert b\vert +\vert c\vert = 2n + 2sN.
$$
If $c_1 TM$ is representable by a K\"ahler $2$-form
(as in the case $M=\C P^n$ or $\grkcn$, for example)  then 
the degree of any holomorphic map is nonnegative,
so it suffices to consider $s=0,1,2,\dots$.  Since 
$N>0$, only a finite number of triple products $\ABC_s$
can be nonzero, so the series defining $\circ$ and $\circ_t$ are
indeed finite series.  
Let us examine briefly the cases $s=0,1,2$.  For $s=0$, we have
\footnote{The following three diagrams, and the one in
\S 4, can be downloaded from http://www.comp.metro-u.ac.jp/\~{}martin}
$\HolD^{\tilde A,p} = \tilde A$, so the definition of $\ABC_s$
reduces to the definition of $\ABC_0$, as it should.  We are simply
counting the points of the triple intersection $\tilde A\cap \tilde B \cap
\tilde C$.

$${}$$
$${}$$
$${}$$
$${}$$

For $s=1$, the triple product $\ABC_1$ counts the 
holomorphic maps $f$ of degree $1$ such that $f(p)\in \tilde A$,
$f(q)\in \tilde B$, $f(r)\in \tilde C$.
In this case,  $\vert a\vert +\vert b\vert +\vert c\vert = 2n + 2N$,
so $\vert a\vert +\vert b\vert +\vert c\vert > 2n$, and the triple
intersection $\tilde A\cap \tilde B \cap \tilde C$ is \ll in general\rr
empty.

$${}$$
$${}$$
$${}$$
$${}$$

For $s=2$,  the triple product $\ABC_2$
counts the holomorphic maps $f$ of degree $2$ such that $f(p)\in \tilde A$,
$f(q)\in \tilde B$, $f(r)\in \tilde C$.

$${}$$
$${}$$
$${}$$
$${}$$

\proclaim{Example: $M=\C P^n$}
\endproclaim

In this case $N=n+1$, so the numerical condition (for $\ABC_s\ne 0$)
is $\vert a\vert +\vert b\vert +\vert c\vert = 2n + 2s(n+1)$. Since
$0\le \vert a\vert,\vert b\vert,\vert c\vert \le 2n$, it follows
immediately that $s=0$ and $s=1$ are the only relevant values,
i.e. that $\ABC_s=0$ for $s\ne 0,1$.  

For $s=0$ we already know that
$$
\lan X_i\vert X_j\vert X_k\ran_0 =
\cases
1\quad\text{if}\ i+j+k=n\\
0\quad\text{otherwise}
\endcases
$$
This is just the triple product for ordinary cohomology.

For $s=1$ we claim that
$$
\lan X_i\vert X_j\vert X_k\ran_1 =
\cases
1\quad\text{if}\ i+j+k=2n+1\\
0\quad\text{otherwise}
\endcases
$$
To prove this, one shows
that, when $i+j+k=2n+1$,
there exist complex linear subspaces $E^i,E^j,E^k$ of $\C^{n+1}$,
of codimensions $i,j,k$,
with the following property:  there exist unique complex
lines $L',L'',L'''$ such that $L'\sub E^i$,  $L''\sub E^j$, $L'''\sub E^k$
and such that $L',L'',L'''$ span a subspace $E$ of dimension $2$.
The holomorphic map of degree $1$ defined by the inclusion
$\P(E)\sub \C P^{n+1}$  is then the unique point of the triple
intersection 
$\Hol_1^{\P(E^i),p} \cap \Hol_1^{\P(E^j),q} \cap \Hol_1^{\P(E^k),r}$,
and we obtain $\lan X_i\vert X_j\vert X_k\ran_1 = 1$.  (The multiplicity
of the intersection point is $1$ because $\Hol_1^{\P(E^i),p}$ is a
smooth subvariety of the space of all holomorphic maps.)

We shall not give a direct proof of the claim, because --- if 
one makes use of the associativity
of quantum multiplication --- an indirect proof is much easier,
as we shall explain later on. This phenomenon is typical: 
a direct (but tedious) problem of \ll enumerating\rr rational
curves with certain properties can sometimes be replaced by
a simple argument making use of the (nontrivial) properties
of quantum cohomology.  Spectacular results of this type have been
obtained by using more general Gromov-Witten invariants, namely the
ones mentioned in the footnote at the beginning of \S 7.

A simple heuristic argument for the calculation of
$\lan X_i\vert X_j\vert X_k\ran_1 $
can be made as follows.  Let $V$ be the linear subspace
of $\C^{n+1}$ of codimension $i$ which is given by the
conditions $z_0=\dots=z_{i-1}=0$. The space
$$
\Hol_s^{\P(V),p} = \{ f=[p_0;\dots;p_n] \st [f]=s\text{ and }
p_0(p)=\dots = p_{i-1}(p)=0 \}
$$
is a dense open subset of a linear subspace
of $\P(\C^{s+1}\times \dots \times \C^{s+1})\cong \P(\C^{(n+1)(s+1)})$ 
of codimension $i$.  (The component functions $p_i$ are complex
polynomials of degree at most $s$, such that no point of 
$\C P^1=\C\cup\infty$ is a common root of $p_0,\dots,p_n$.)
In the computation of 
$\lan X_i\vert X_j\vert X_k\ran_1$ we are therefore considering
the intersection of (dense open subsets of) three linear
subspaces of the $(2n+1)$-dimensional projective
space $\P(\C^{2(n+1)})$ of codimensions $i,j,k$. We expect
a single intersection point if and only if $i+j+k=2n+1$,
at least if the three subspaces are in general position and if the
intersection point can be assumed to lie in each of the respective
open sets. Thus the computation of  $\lan X_i\vert X_j\vert X_k\ran_1$
is, like the computation of  $\lan X_i\vert X_j\vert X_k\ran_0$,
a problem of finding intersections of linear subspaces --- the 
difficulty being the fact that we are dealing with dense open subspaces,
rather than the linear subspaces themselves.

Let us calculate the quantum products $x_i\circ x_j$.  We have
$$
x_i\circ x_j = (x_i\circ x_j)_0 + (x_i\circ x_j)_1q
$$
where
$\vert  (x_i\circ x_j)_0\vert = 2i+2j$, 
$\vert  (x_i\circ x_j)_1\vert = 2i+2j-2(n+1)$, and
$$
\align
\lan (x_i\circ x_j)_0, X_k\ran &= 
\lan X_i \vert X_j \vert X_k \ran_0 =
\cases
1\quad\text{if}\ i+j+k=n\\
0\quad\text{otherwise}
\endcases
\\
\lan (x_i\circ x_j)_1, X_k\ran &= 
\lan X_i \vert X_j \vert X_k \ran_1 =
\cases
1\quad\text{if}\ i+j+k=2n+1\\
0\quad\text{otherwise}
\endcases
\endalign
$$
It follows that
$$
\align
(x_i\circ x_j)_0 &= x_{i+j}\quad\text{if}\ 0\le i+j\le n\\
(x_i\circ x_j)_1 &= x_{i+j-(n+1)}\quad\text{if}\ n+1\le i+j\le 2n
\endalign
$$
and we conclude that
$$
x_i\circ x_j =
\cases
x_{i+j}\quad\text{if}\ 0\le i+j\le n\\
x_{i+j-(n+1)}q\quad\text{if}\ n+1\le i+j\le 2n
\endcases
$$

The (small) quantum cohomology algebra is therefore
$$
QH^\ast(\C P^n;\C)\ \cong\ 
\C[p,q,q^{-1}]/\lan qq^{-1}-1,\ p^{n+1}-q\ran
$$
where we have written $p$ for $x_1$, and $1$ for $x_0$.  The subalgebra
$$
\widetilde{QH}^\ast(\C P^n;\C)\ \cong\ 
\C[p,q]/\lan  p^{n+1}-q\ran
$$ 
is obtained by considering only nonnegative powers of $q$.

Notice that $QH^\ast(\C P^n;\C)$, in contrast to $(H^\ast(\C P^n;\C),\circ_t)$,
is an infinite-dimensional vector space. However, 
$QH^\ast(\C P^n;\C)$ contains the same information as $(H^\ast(\C P^n;\C),\circ_t)$,  
because of the \ll periodicity isomorphisms\rr
$$
\times q: QH^i(\C P^n;\C) \to QH^{i+2n+2}(\C P^n;\C).
$$
To illustrate this, we list additive generators of
$QH^i(\C P^n;\C)$ for various
degrees $i$ in the table below:
$$
{\eightpoint
\matrix
\dots & \ssize{-2n-2} & \ssize{-2n} & \dots & \ssize{-2} & \ssize{0} & \ssize{2} & 
\dots & \ssize{2n} & \ssize{2n+2} & \ssize{2n+4} & \dots & \ssize{4n+2} & \dots\\
\dots  & x_0q^{-1} & x_1q^{-1} & \dots & x_nq^{-1} & x_0 & x_1 & 
\dots & x_n & x_0q & x_1q & \dots & x_nq &\dots
\endmatrix}
$$
Notice also that the \ll obvious\rr product structure on 
$H^\ast(\C P^n;\C)\otimes \La$ results in the algebra
$\C[p]/\lan p^{n+1}\ran \otimes \La$; this is quite different
from the algebra $\C[p]\otimes\La/\lan p^{n+1} - q\ran$ obtained from the
quantum product. 

We can now explain the (deceptively trivial) indirect calculation of 
the quantum cohomology of $\C P^n$, which is the one usually given in 
expositions of the subject. The point is that, for dimensional reasons, 
\ll the  only nontrivial quantum product is $x_1\circ x_n = q$\rrr. 
For example, since $\vert q\vert= 2n+2$, there can be no term involving
$q$ in the quantum product $x_i\circ x_j$ when $i+j<n$, so we must have
$x_i\circ x_j=x_{i+j}$. From this and the formula $x_1\circ x_n = q$,
all other quantum products $x_i\circ x_j$ may be deduced, e.g.
$x_2\circ x_n = x_1^2\circ x_n = (x_1\circ x_1)\circ x_n = 
x_1\circ (x_1\circ x_n) = x_1q$.
To establish the formula $x_1\circ x_n = q$,
one must prove that $\lan X_1\vert X_n\vert X_n\ran_1=1$. This is easy,
as $\lan X_1\vert X_n\vert X_n\ran_1$ is the number of linear maps
$\C P^1 \to \C P^n$ which \ll hit\rr generic representatives of
$X_1$ (a hyperplane), $X_n$ (represented by a point not on the hyperplane), 
$X_n$ (represented by another point not on the hyperplane) at three
prescribed points of $\C P^1$.  But this argument is valid only if we assume
various properties of the quantum product, such as its associativity and
the fact that it is a deformation of the cup product.  On the other hand, if we
calculate all Gromov-Witten invariants $\lan X_i\vert X_j\vert X_k\ran_s$
directly by linear algebra, then we will (in a very inefficient manner)
establish these properties of the quantum product for the case of $\C P^n$. 

We conclude with some comments on the technical conditions which are
needed to justify the naive definition of $\ABC\subD$.  First, 
a very general definition of $\ABC\subD$ (and hence of the quantum
product)  is possible under the assumption that the (connected simply
connected K\"ahler---or even merely symplectic) manifold $M$ is 
\ll positive\rr in some sense, for example that $\lan c_1 TM,D\ran > 0$ 
for each homotopy class $D\in \pi_2(M)$
which contains a holomorphic map $\C P^1 \to M$.  A Fano manifold --- that is,
a manifold $M$ for which the cohomology class $c_1 TM$ can be 
represented by a K\"ahler $2$-form --- is automatically positive in this sense.
It can be shown that the quantum product is commutative and associative, 
when $M$ is positive.  So far, however, there is no guarantee that the
Gromov-Witten invariant $\ABC\subD$ can be {\it computed} by the naive
formula given earlier in this section. For this, one needs an additional
assumption, for example that $M$ is \ll convex\rr in the sense that
$H^1(\C P^1,f^\ast TM)=0$ for all holomorphic  $f:\C P^1\to M$. Convexity
implies in particular that $\HolD$ is a manifold, and that it has the 
\ll expected\rr  dimension $n + \lan c_1 TM, D\ran$. Homogeneous K\"ahler
manifolds are convex, for example.

The series defining the quantum product will in general contain
infinitely many powers of $q$ (both positive and negative).  However, if
there exist $D_1,\dots,D_r\in \pi_2(M)$ such that all holomorphically
representable classes $D$ are of the form $\sum_{i=1}^r n_i D_i$
with $n_i\ge 0$, then it follows from the positivity condition that each
such series contains only a finite number of terms $q^D=q_1^{s_1}\dots q_r^{s_r}$,
and that $s_i\ge 0$ in each case. Any simply connected homogeneous
K\"ahler manifold satisfies this condition, as well as positivity and
convexity.  However, problems arise as soon as we contemplate
nonhomogeneous manifolds --- even for such simple examples as the
Hirzebruch surfaces (see Appendix 2).

\head
\S 4 Landau-Ginzburg potentials
\endhead

References: \cite{Si-Ti}, \cite{Be1}-\cite{Be3}, \cite{RRW}

There is a rather surprising \ll analytic\rr description of the
cohomology algebra of the Grassmannian $\grkcn$ (and certain other Hermitian symmetric
spaces), which was discovered recently by physicists (\cite{LVW}).
Before stating this, we shall review the two standard
descriptions, which are well known to topologists.  We shall
generally follow the notation of \cite{Si-Ti}.

The \ll algebraic\rr (Borel) description of this cohomology algebra is
$$
H^\ast(\grkcn;\Z) \cong
\frac{\Z[u_1,\dots,u_k,v_1,\dots,v_{n-k}]^{\Sigma_{k,n-k}}}
{\Z[u_1,\dots,u_k,v_1,\dots,v_{n-k}]^{\Sigma_n}_+}
$$
where  
$\Z[u_1,\dots,u_k,v_1,\dots,v_{n-k}]^{\Sigma_{k,n-k}}$ denotes the
polynomials which are (separately) symmetric
in $u_1,\dots,u_k$ and in $v_1,\dots,v_{n-k}$, 
and $\Z[u_1,\dots,u_k,v_1,\dots,v_{n-k}]^{\Sigma_n}_+$
denotes the  polynomials of positive degree which are
symmetric in $u_1,\dots,u_k,v_1,\dots,v_{n-k}$.  Each
of $u_1,\dots,u_k,v_1,\dots,v_{n-k}$ represents
a cohomology class of degree $2$ here.
This description is equivalent to the slightly more geometrical
description
$$
H^\ast(\grkcn;\Z) \cong
\Z[c_1,\dots,c_k,s_1,\dots,s_{n-k}]/
\lan c(\Cal V) c(\C^n/\Cal V) - 1 \ran
$$
where
$$
\align
c(\Cal V) &= 1+c_1+c_2+\dots+c_k \\
c(\C^n/\Cal V) &= 1+s_1+s_2+\dots+s_{n-k}
\endalign
$$
are the total Chern classes of the tautologous bundle $\Cal V$ 
(of rank $k$) and the quotient bundle $\C^n/\Cal V$ (of rank $n-k$)
on $\grkcn$.   The previous description is obtained by writing
$1+c_1+c_2+\dots+c_k = (1+u_1)\dots(1+u_k)$ and
$1+s_1+s_2+\dots+s_{n-k} = (1+v_1)\dots(1+v_{n-k})$,
in accordance with the Splitting Principle 
(see chapter 4 of \cite{Bo-Tu}).

Taking the components in each positive degree of the identity 
$c(\Cal V) c(\C^n/\Cal V) = 1$, we obtain $n$ relations between
$c_1,\dots,c_k$ and $s_1,\dots,s_{n-k}$.  To express these
relations, it is convenient to consider them as the first $n$
in an infinite series of formal relations 
$$
(1+c_1+c_2+\dots)(1+s_1+s_2+\dots) = 1
$$
obtained from the product of two formal power series. 
We can use the first $n-k$ relations to express $s_1,\dots,s_{n-k}$
in terms of $c_1,\dots,c_k$; the next $k$ relations can
then be written in the form
$$
s_{n-k+1}=f_{n-k+1}(c_1,\dots,c_k),\quad\dots,\quad s_{n}=f_{n}(c_1,\dots,c_k)
$$
for some polynomials $f_{n-k+1},\dots,f_{n}$.  It follows that
$$
H^\ast(\grkcn;\Z) \cong
\Z[c_1,\dots,c_k]/
\lan  f_{n-k+1},\dots,f_{n}    \ran.
$$
Alternatively, we can use the first $k$ relations to express $c_1,\dots,c_{k}$
in terms of $s_1,\dots,s_{n-k}$; the next $n-k$ relations are then of the form
$$
c_{k+1}=g_{k+1}(s_1,\dots,s_{n-k}),\quad\dots,\quad c_{n}=g_{n}(s_1,\dots,s_{n-k})
$$
for some polynomials $g_{k+1},\dots,g_{n}$. We obtain
$$
H^\ast(\grkcn;\Z) \cong
\Z[s_1,\dots,s_{n-k}]/
\lan  g_{k+1},\dots,g_{n}    \ran.
$$

\proclaim{Example: $M=\C P^n$}
\endproclaim

Here $k=1$ and we have $f_{n+1}(c_1)=(-1)^{n+1} c_1^{n+1}$.  This gives
the familiar description
$$
H^\ast(\C P^n;\Z) \cong
\Z[c_1]/
\lan   c_1^{n+1}   \ran.
$$
Alternatively, in terms of $s_1,\dots,s_n$, we have a rather more
complicated description
$$
H^\ast(\C P^n;\Z) \cong
\Z[s_1,\dots,s_{n}]/
\lan g_{2},\dots,g_{n+1}   \ran.
$$
The relations $g_{2}=0,\dots,g_{n+1}=0$ are equivalent to (but not identical
to)  $s_{i+1}=s_1s_i$, $1\le i\le n$, which of course agrees with the
previous description.

\proclaim{Example: $M=Gr_2(\C^4)$}
\endproclaim

Here it makes no difference whether we use $c_1,c_2$ or $s_1,s_2$.
Let us choose $c_1,c_2$.  Then the relations are
$$
\align
s_1&= -c_1\\
s_2&= -s_1c_1 - c_2 = c_1^2 - c_2\\
s_3&= -s_2c_1 -s_1c_2 - c_3= -c_1^3 + 2c_1c_2 -c_3\\
s_4&= -s_3c_1 -s_2c_2 -s_1c_3 - c_4 = c_1^4 -3c_1^2c_2 + c_2^2 + c_1c_3 - c_4
\endalign
$$
Setting $c_3=c_4=0$, we obtain
$$
f_3(c_1,c_2)=-c_1^3 + 2c_1c_2,\quad f_4(c_1,c_2)= c_1^4 -3c_1^2c_2 + c_2^2.
$$
The cohomology algebra is therefore
$$
H^\ast(Gr_2(\C^4);\Z) \cong
\Z[c_1,c_2]/
\lan f_{3},f_{4}   \ran.
$$
Notice that $f_4+c_1 f_3 = c_2^2 - c_1^2c_2$, so we can replace $f_4$
by this to obtain the slightly simpler description
$$
H^\ast(Gr_2(\C^4);\Z) \cong
\Z[c_1,c_2]/
\lan -c_1^3 + 2c_1c_2, c_2^2 - c_1^2c_2   \ran.
$$
This could have been obtained immediately from the identity
$(1+c_1+c_2+\dots)(1+s_1+s_2+\dots) = 1$ by setting
$c_i=s_i=0$ for $i>2$. However, we prefer not to do this,
as the polynomials $f_i$ and $g_j$ are more convenient for the
general theory.

Now we turn to the \ll geometrical\rr (Schubert)
description of the cohomology algebra of $\grkcn$.  The rank of
the abelian group $H^\ast(\grkcn;\Z)$ is $\binom nk$.
Additive generators may be parametrized by Young diagrams
of the form
$${}$$
$${}$$
$${}$$
$${}$$
$${}$$
$${}$$
or (less picturesquely) by $k$-tuples $(\la_1,\dots,\la_k)\in\Z^k$
such that $0\le \la_k \le \la_{k-1}\le \dots \le \la_{1} \le n-k$.
Observe that there are $\binom nk$ such diagrams, because they are
in one to one correspondence with choices of $k$ distinct elements
$\la_k +1,\la_{k-1}+2,\dots,\la_{1}+k$ from the set $\{1,2,\dots,n\}$.
The Schubert variety $X(\la)$ is defined to be the subvariety
$$
X(\la)=\{ V\in\grkcn \st \dim V\cap \C^{n-k+i-\la_i}\ge i, 1\le i\le k \}
$$
of $\grkcn$.  It is an irreducible subvariety of complex codimension
$\sum_{i=1}^k \la_i$, and therefore it defines a homology class.  The
Poincar\acuteaccent e dual cohomology class, of degree 
$2\sum_{i=1}^k \la_i$, will be denoted by $x(\la)$.  

Giambelli's formula expresses $x(\la)$ in terms of the multiplicative
generators $s_1,\dots,s_{n-k}$ of the Borel description:
$$
x(\la)=
\left |
\matrix
s_{\la_{1}}  &  s_{\la_{1}+1}  &  s_{\la_{1}+2}  & \dots & s_{\la_{1}+k-1} \\
s_{\la_{2}-1}  &  s_{\la_{2}}  &  s_{\la_{2}+1}  & \dots & s_{\la_{2}+k-2} \\
\dots  &  \dots  &  \dots  & \dots & \dots \\
\dots  &  \dots  &  \dots  & \dots & \dots \\
s_{\la_{k}-k+1}  &  s_{\la_{k}-k+2}  &  s_{\la_{k}-k+3}  & \dots & s_{\la_{k}} \\
\endmatrix
\right |
$$
(where $s_i$ is defined to be zero if $i\notin\{0,1,\dots,n-k\}$,
and $s_0=1$).  These polynomials in $s_1,\dots,s_{n-k}$ (parametrized
by $\la$) are known as Schur functions.
As a particular case of this formula we have
$$
s_j = x(j,0,\dots,0), \quad 1\le j\le n-k
$$
and the Poincar\acuteaccent e dual Schubert variety is
$$
X(j,0,\dots,0) = \{ V\in\grkcn \st \dim V\cap \C^{n-k+1-j} \ge 1 \}.
$$
Such varieties are known classically as the 
\ll special\rr Schubert varieties. We also have
$$
c_j = (-1)^jx(1,\dots,1,0,\dots,0) 
\text{ (with $j$ entries equal to $1$)}, \quad 1\le j\le k
$$
(see \cite{Si-Ti}).

Intersection of Schubert varieties (after translation by
generic elements of $\glnc$) corresponds to multiplication of
the corresponding cohomology classes --- this is the Schubert calculus.
In addition to the Giambelli formula, there is an explicit combinatorial
formula for the product $x(\la)x(\la^\prime)$, called the
Littlewood-Richardson rule.  For the special products
$x(j,0,\dots,0)x(\la^\prime)$ this is a classical formula, called
Pieri's formula.  Giambelli's formula and Pieri's formula together
determine all the products $x(\la)x(\la^\prime)$.

The analytic (Landau-Ginzburg) description of the cohomology algebra
is a re-interpretation of the Borel description.  It is
$$
H^\ast(\grkcn;\Z) \cong
\Z[c_1,\dots,c_k]/
\lan  dP    \ran,
$$
where $P$ is a certain $\C$-valued polynomial function of 
$c_1,\dots,c_k$.  In other words, it is asserted that the
$\C^k$-valued function $(f_{n-k+1},\dots,f_{n})$ has
a \ll primitive\rrr.  Although this might seem superficial, it
turns out that the existence of such \ll Landau-Ginzburg potentials\rr
is fundamental in homological geometry.

\proclaim{Example: $M=\C P^n$}
\endproclaim

A Landau-Ginzburg potential here is 
$P(c_1)= \pm c_1^{n+2}/(n+2)$.

\proclaim{Example: $M=Gr_2(\C^4)$}
\endproclaim

A Landau-Ginzburg potential here is
$$
P(c_1,c_2)= \frac15( c_1^5 - 5c_1^3 c_2 + 5 c_1 c_2^2 ), 
$$
as one verifies by differentiation:
$$
\align
\frac{\bd P}{\bd c_1} &= c_1^4 - 3c_1^2c_2 + c_2^2 = f_4 \\
\frac{\bd P}{\bd c_2} &= -c_1^3 + 2c_1c_2 = f_3.
\endalign
$$

The function $P$ first arose in physics, and
from a mathematical point of view it is quite
mysterious.  However --- with hindsight ---  it is easy
to establish the existence of $P$, as follows.
Consider the formal power series
$$
\align
C(t) &= 1 + c_1 t + c_2 t^2 + \dots \\
S(t) &= 1 + s_1 t + s_2 t^2 + \dots
\endalign
$$
with $C(t)S(t)=1$.  Let
$$
\log C(t) =W(t)= w_1 t + w_2 t^2 + \dots .
$$
We have
$$
\frac{\bd}{\bd c_i} W(t) =
C(t)^{-1} \frac{\bd}{\bd c_i} C(t) =
S(t) t^i.
$$
Equating the coefficients of $t^{n+1}$ here, we obtain
$$
\frac{\bd w_{n+1}}{\bd c_i} = s_{n+1-i}, \quad (1\le i\le k).
$$
Interpreting this in the Borel description of the cohomology algebra,
we see that a  Landau-Ginzburg potential is just
$$
P=w_{n+1}.
$$
We can express this more explicitly in terms of the original
variables $u_1,\dots,u_k$.  Recall that these were defined
by $1+c_1+c_2+\dots+c_k = (1+u_1)\dots(1+u_k)$.  We have
$$
\align
w_1 t + w_2 t^2 + \dots &=
\log(1 + c_1 t + c_2 t^2 + \dots +c_kt^k) \\
&= \sum_{i=1}^k \log(1+u_it) \\
&= \sum_{i=1}^k (-u_i t + \frac12 u_i^2 t^2 - \dots ),
\endalign
$$
hence
$$
P=w_{n+1}= \frac{(-1)^{n+1}}{n+1} \sum_{i=1}^k u_i^{n+1}.
$$

Now we turn to quantum cohomology.  An elementary calculation
(making use of the properties of the quantum product --- see \cite{Si-Ti}) gives:
$$
\widetilde{QH}^\ast(\grkcn;\Z) \cong
\Z[c_1,\dots,c_k,q]/
\lan  f_{n-k+1},\dots,f_{n-1},f_{n} + (-1)^{n-k} q    \ran.
$$

Hence --- by inspection --- we have
$$
\widetilde{QH}^\ast(\grkcn;\Z) \cong
\Z[c_1,\dots,c_k,q]/
\lan  d\tilde P  \ran,
$$
where
$$
\tilde P  = P + (-1)^{n-k} c_1 q = 
\frac{(-1)^{n+1}}{n+1} \sum_{i=1}^k u_i^{n+1} + (-1)^{n-k}\sum_{i=1}^k u_i.
$$
Thus, we have a Landau-Ginzburg potential for quantum cohomology too.

A nontrivial application of  Landau-Ginzburg potentials is the
existence of an analytic formula for quantum products (and ordinary
products).  This can be expressed either as an integral or a
sum of residues.  Moreover,
the formula generalizes to a (conjectural) formula for higher genus
Gromov-Witten invariants, the \ll Formula of Vafa and Intriligator\rrr.

One version of the formula is as follows (see \cite{Be1}).  Let 
$T(c_1,\dots,c_k)$ be a polynomial whose total degree is equal
to the dimension of $\grkcn$.  This is an integer multiple of
the generating class, and the integer is given up to sign by
$$
 \pm \sum_{x} \frac {T(x)}{h(x)}
$$
where the sum is over the (finite) set of values $x$ such
that $d \tilde P(x)=0$, and where
$$
h = \det\left( \frac{ \bd^2 \tilde P}{\bd c_i \bd c_j} \right).
$$

\head
\S 5 Homological geometry
\endhead

References: \cite{Gi-Ki}, \cite{Au2}-\cite{Au4}

We shall consider a sequence of spaces, denoted (i)-(iv) below.
As we progress, we shall make the transition from the algebra
of homology/cohomology theory to the geometry of manifolds/varieties.
In this way we enter the world of \ll homological geometry\rrr,
a world (see section 1 of \cite{Gi1}) where one deals with 
functions (and differential geometric
objects) defined on homology and cohomology vector spaces.

\no(i) The manifold $M$.

We assume that $M$ satisfies the technical conditions of \S 3.
To proceed further, we shall need the cohomology 
$R$-modules $H^\ast(M;R)$ for
$R=\Z$,
$R=\C$, and 
$R=\C^\ast=\C-\{0\}\cong\C/2\pi\i\Z$. 

\no(ii) The manifold $B$.

Let us fix additive bases as follows:
$$
\align
H_2(M;\Z) &= \bigoplus_{i=1}^r \Z A_i\\
H^2(M;\Z) &= \bigoplus_{i=1}^r \Z b_i
\endalign
$$ 
A typical element of $H_2(M;\C)$ will be denoted by
$\sum_{i=1}^r p_i A_i$.  Thus, we regard $p_i$ as the $i$-th
\ll coordinate function\rr on  $H_2(M;\C)$. 
Similarly, a typical element of $H^2(M;\C)$ will be denoted by
$\sum_{i=1}^r t_i b_i$.

We introduce the complex algebraic torus
$$
B= H^2(M;\C/2\pi\i \Z) = \bigoplus_{i=1}^r \C^\ast [b_i] 
\cong \C^\ast \times \dots \times \C^\ast
$$
where $[b_i]$ is the element of 
$H^2(M;\C/2\pi\i \Z) \cong H^2(M;\C) / H^2(M;2\pi\i \Z)$
corresponding to $b_i$.
A typical element of $B$ will be denoted by $(q_1,\dots,q_r)$;
thus we regard $q_i$ as the $i$-th coordinate function on $B$.

Unfortunately the notation $p_i,t_i,q_i$ here
breaks with our tradition of using lower-case letters for 
cohomology classes. In the case of $p_i$ there is really no
problem as it is a linear functional on $H_2(M;\C)$, and therefore can be
identified with an element of $H^2(M;\C)$. To make this more
concrete, let us choose $b_i=a^c_i$, so that
$\lan b_i,A_j\ran = \lan a^c_i,A_j\ran = \lan a^c_ia_j,M\ran = \de_{ij}$.
Then $\lan b_i,\ \ran$ is just the $i$-th coordinate function
$p_i$.  In other words, the cohomology class $b_i\in H^2(M;\C)$
corresponds to the coordinate function $p_i\in  H_2(M;\C)^\ast$.

In the case of $q_i$, which previously denoted a \ll formal variable\rrr, we
can justify the current functional interpretation as follows.  
Recall (from \S 2) that the two versions of quantum cohomology are
related by \ll putting $q^D=e^{\lan t,D\ran}$\rrr, and hence
\ll $q_i=e^{\lan t,D_i\ran}$\rrr. More precisely, in the first version
(using $\circ$) we regard $q$ as a formal variable, whereas in the second
version (using $\circ_t$) we regard $q_i$ as the function 
$t\mapsto e^{\lan t,D_i\ran}$ on $H^2(M;\C)$.
Let us choose the basis $A_1,\dots,A_r$ to be of the type denoted
previously (at the ends of \S 2 and \S 3) by $D_1,\dots,D_r$. Then
our current notation $q_i$ denotes the function
$$
q_i:H^2(M;\C)/ H^2(M;2\pi\i \Z) \to \C^\ast,
\quad
[t]\mapsto e^{ \lan t,A_i\ran},
$$
which is (induced by) our earlier realization
of the formal variable $q_i$ as a function.

In the case of the notation $\sum_{i=1}^r t_i b_i$, we shall just
have to ask the reader to bear in mind that $t_i$ is a coordinate function on
$H^2(M;\C)$, not a cohomology class.  This has the advantage, at least,
that we can write $t=\sum_{i=1}^r t_i b_i$ for a general cohomology
class in $H^2(M;\C)$. (In contrast, we must be careful
to avoid writing $p=\sum_{i=1}^r p_i A_i$ or $p=(p_1,\dots,p_r)$, 
because $(p_1,\dots,p_r)$ is a homology class, and $p$ is
a cohomology class. Similarly, we shall
avoid writing  $q=(q_1,\dots,q_r)$, because $(q_1,\dots,q_r)$
is an element of $B$. The case $r=1$ will need special care!)

These remarks on notation may seem tiresome, but we want to emphasize that
we shall regard $p_i$ and $q_i$ primarily as {\it coordinate functions} on
$H_2(M;\C)$ and $B$ (respectively), from now on.  As is usual, 
we shall then write $\bd/\bd p_i$, $\bd/\bd q_i$
for the vector fields associated to these coordinate functions.

\no(iii) The manifold $T^\ast B$.

As $B$ is a group, we have canonical isomorphisms
$$
\align
TB&\cong B \times H^2(M;\C)\\
T^\ast B &\cong B \times H^2(M;\C)^\ast.
\endalign
$$
Like all cotangent bundles, $T^\ast B$ has a natural symplectic 
form $d\la$, where $\la$ is the Liouville form.  The symplectic
manifold $T^\ast B$ will be the focus of
our attention in this section.  We  have
$$
T^\ast B \cong B \times H_2(M;\C)
$$
if we use the above identification
$$
H^2(M;\C) \cong H_2(M;\C)^\ast, \quad
b_i\mapsto p_i=\lan b_i, \ \ran.
$$
Then $p_i$ and $q_i$ may be regarded as $\C$-valued functions on
the manifold $T^\ast B$.
With this terminology, the $1$-form $\la$ is given explicitly by
$$
\la = \sum_{i=1}^r p_i\wedge \frac{dq_i}{q_i}.
$$

Next we introduce the group algebra
$$
\align
\La &= \C[H_2(M;\Z)] \\
&= \{ \text{polynomials in $q_1,\dots,q_r,q_1^{-1},\dots,q_r^{-1}$} \} \\
&=\{ \text{regular functions on $B$} \}
\endalign
$$
and the symmetric algebra
$$
\align
S &= S(H^2(M;\C)) \\
&\cong S(H_2(M;\C)^\ast) \\
&=\{ \text{polynomials in $p_1,\dots.p_r$} \}\\
&=\{ \text{regular functions on $H_2(M;\C)$} \}.
\endalign
$$
It follows that the algebra $S\otimes \La$ may be identified
with the algebra of regular functions on $T^\ast B$.

For the rest of this section we shall {\it assume} that 
\ll $H^2(M;\Z)$ generates $H^\ast(M;\Z)$\rrr, 
so that the natural homomorphism
$$
S \to H^\ast(M;\C)
$$
is surjective.  Hence $H^\ast(M;\C)$ has the form
$$
H^\ast(M;\C)\cong S/I
\cong \C[p_1,\dots,p_r]/\lan R_1,R_2,\dots \ran
$$
where the ideal $I$ (generated by $R_1,R_2,\dots$)
is the kernel of the above homomorphism.  
The manifolds $\C P^n$ and $F_n=F_{1,2,\dots,n-1}(\C^n)$ satisfy
this assumption, but $\grkcn$ (for $2\le k \le n-2$) does not.

It follows that the natural homomorphism
$$
S\otimes \La \to QH^\ast(M;\C)
$$
is surjective as well, and hence (by Theorem 2.2 of \cite{Si-Ti})
$$
QH^\ast(M;\C)\cong S\otimes \La/ \Cal I
\cong \C[p_1,\dots,p_r,q_1,q_1^{-1},\dots,q_r,q_r^{-1}]/
\lan \Cal R_1,\Cal R_2,\dots \ran
$$
for some relations $\Cal R_1,\Cal R_2,\dots$ which are \ll quantum
versions\rr of $R_1,R_2,\dots$ .

\no(iv) The variety $V_M$.

Since $S\otimes \La$ is the \ll coordinate ring\rr of $T^\ast B$,
the quotient ring $S\otimes \La / \Cal I$ should be the
coordinate ring of a subvariety $V_M$ of $T^\ast B$.  More
informally, we shall just
consider $V_M$ to be the subvariety of $T^\ast B$
defined by the equations $\Cal R_i=0$, i.e.
$$
V_M=\{ (q_1,\dots,q_r,p_1,\dots,p_r) \in T^\ast B
\st \Cal R_1 = \Cal R_2 = \dots = 0 \}
$$

It is shown in \cite{Gi-Ki} and \cite{Au4}
that, under certain conditions,  $V_M$
is a Lagrangian subvariety of $T^\ast B$.  This means that $V_M$
is maximal isotropic with respect to the symplectic form $d\la$.
We shall discuss this result in more detail in \S 7.  For
the moment, we consider two related properties:

\no(L1)\quad  $\Cal R,\Cal S \in \Cal I 
\implies \{\Cal R,\Cal S\}\in\Cal I$

\no(L2)\quad  $\Cal R,\Cal S \in \Cal I 
\implies \{\Cal R,\Cal S\}=0$

\no where $\{\ ,\ \}$ is the Poisson bracket associated to
the symplectic form $d\la$.  

To explain the relevance of these properties, we need to review
some facts on symplectic geometry (see \cite{We}). First, any
(smooth) function $f:T^\ast B\to \C$ has an associated
\ll Hamiltonian\rr vector field $H_f$, defined by
$H_f = (d\la)^{-1}\circ df$. The Poisson bracket $\{f,g\}$ of two
such functions is defined by $\{f,g\}=d\la(H_g,H_f)$
(this is equal to $dg(H_f)$ or $-df(H_g)$). In terms of
the coordinates $p_1,\dots,p_r,q_1,\dots,q_r$ we have
$$
\{f,g\} = -\sum_{i=1}^r q_i
\left(
\frac{\bd f}{\bd p_i}\frac{\bd g}{\bd q_i} -
\frac{\bd g}{\bd p_i}\frac{\bd f}{\bd q_i}.
\right)
$$

We shall assume that (the smooth part of) $V_M$ is 
a regular level set of the functions in $\Cal I$, i.e. that
$$
T_m V_M = \bigcap_{\Cal S\in\Cal I} \Ker d\Cal S_m
$$
for all smooth points $m$. 

Now, let us assume in addition that property (L1) holds.  Then
we have
$$
d\Cal S(H_{\Cal R})\vert_{V_M} = \{\Cal R,\Cal S\}\vert_{V_M} = 0
$$
for all $\Cal R,\Cal S\in\Cal I$, so each vector field $H_{\Cal R}$
is tangent to $V_M$.  From the well known formula
$[H_g,H_f]=H_{\{f,g\}}$, and property (L1) again, it follows that
$$
[H_{\Cal R}\vert_{V_M},H_{\Cal S}\vert_{V_M}]=
[H_{\Cal R},H_{\Cal S}]\vert_{V_M}=
H_{\{\Cal S,\Cal R\}}\vert_{V_M}=
0
$$
for all $\Cal R,\Cal S\in\Cal I$.  Hence, the vector fields
$H_{\Cal R}\vert_{V_M}$ define an integrable distribution --- with
\ll linear\rr leaves --- on the smooth part of $V_M$.  Any 
integral manifold $V$ is isotropic with respect to $d\la$, since
$d\la(H_{\Cal R},H_{\Cal S})\vert_{V}=\{\Cal S,\Cal R\}\vert_V = 0$.

In general, $V$ is strictly smaller than $V_M$; in fact if
$\Cal I$ is generated by relations $\Cal R_1,\dots,\Cal R_{r^\pr}$,
then we have $\dim V\le {r^\pr}$ and $\codim V_M\le {r^\pr}$.  However, in
the special case where $V=V_M$ (and in particular $r={r^\pr}=\dim B$),
then our discussion shows that property (L1) implies that
$V_M$ is isotropic and in fact maximal isotropic, i.e. Lagrangian.

As a temporary substitute for the general proof that $V_M$ is Lagrangian,
we shall verify property (L1) for several examples.  In addition, we
shall consider whether or not property (L2) is satisfied. 

\proclaim{Example: $M=\C P^n$}
\endproclaim

We have $V_M=\{(q,p)\in \C^\ast\times\C \st p^{n+1}=q \}$.  Thus
$V_M$ is a smooth submanifold of $\C^\ast\times\C$ (it is isomorphic
to $\C^\ast$).  Being one-dimensional, it is automatically Lagrangian.
Property (L1), and indeed the stronger property (L2), is automatically
satisfied in this example.

\proclaim{Example: $M=F_{1,2}(\C^3)$}
\endproclaim

From  Appendix 1, $V_M$ is the subvariety of
$$
T^\ast B=
\{(q_1,q_2,x_1,x_2,x_3)\in (\C^\ast)^2\times \C^3 \st x_1+x_2+x_3=0\}
\cong (\C^\ast)^2\times \C^2
$$
defined by the equations
$$
\align
\Cal R_1 &=  x_1x_2+x_2x_3+x_3x_1+q_1+q_2 = 0\\
\Cal R_2 &=  x_1x_2x_3+x_3q_1+x_1q_2 = 0
\endalign
$$
We shall now make the following change of coordinates:
$$
p_1=x_1,\quad
p_2=x_1+x_2
$$
(this is dictated by our choice of isomorphism
$T^\ast B \cong B\times H_2(M;\C)$, i.e. we require
$\lan p_i,A_j\ran = \de_{ij}$, where the $q_i$'s are defined relative
to the basis consisting of the $A_i$'s).  
Then we obtain
$$
\align
\Cal R_1 &=  -p_1^2 -p_2^2 +p_1p_2 +q_1+q_2\\
\Cal R_2 &=  -p_1p_2^2 +p_1^2p_2 -p_2q_1 +p_1q_2
\endalign
$$
Computing derivatives, we have
$$
\align
\frac{\bd \Cal R_1}{\bd p_1}\frac{\bd \Cal R_2}{\bd q_1} -
\frac{\bd \Cal R_2}{\bd p_1}\frac{\bd \Cal R_1}{\bd q_1}
&= (-2p_1+p_2)(-p_2)-(-p_2^2+2p_1p_2+q_2) = -q_2\\
\frac{\bd \Cal R_1}{\bd p_2}\frac{\bd \Cal R_2}{\bd q_2} -
\frac{\bd \Cal R_2}{\bd p_2}\frac{\bd \Cal R_1}{\bd q_2}
&= (-2p_2+p_1)(p_1)-(-2p_1p_2+p_1^2-q_1) = q_1,
\endalign
$$
so
$$
\{ \Cal R_1,\Cal R_2\} = -q_1(-q_2) - q_2(q_1) = 0.
$$
Thus, the stronger condition (L2) is satisfied in this case.

\proclaim{Example: $M=\Si_1 = \P(\Cal O(0)\oplus \Cal O(-1) )$}
\endproclaim

With the notation of Appendix 2, the quantum cohomology algebra of 
the Hirzebruch surface $\Si_1$ is
$$
\widetilde{QH}^\ast(\Si_1;\C) \cong
\frac{\C[x_1,x_4,q_1,q_2]}
{\lan  x_1^2 - x_2 q_2,  x_4^2 - z - q_1  \ran}
$$
where $x_2=x_4-x_1$ and $z=x_1x_4$.
We introduce the new notation
$$
p_1=x_4,\quad p_2=x_1
$$
(again this is dictated by our choice of isomorphism
$T^\ast B \cong B\times H_2(M;\C)$),
and obtain the relations
$$
\align
\Cal R_1 &=  p_2^2 - (p_1-p_2) q_2\\
\Cal R_2 &= p_1^2 - p_1p_2 - q_1.
\endalign
$$
Computing derivatives, we have
$$
\align
\frac{\bd \Cal R_1}{\bd p_1}\frac{\bd \Cal R_2}{\bd q_1} -
\frac{\bd \Cal R_2}{\bd p_1}\frac{\bd \Cal R_1}{\bd q_1}
&= (-q_2)(-1)-(2p_1-p_2)(0)=q_2\\
\frac{\bd \Cal R_1}{\bd p_2}\frac{\bd \Cal R_2}{\bd q_2} -
\frac{\bd \Cal R_2}{\bd p_2}\frac{\bd \Cal R_1}{\bd q_2}
&= (2p_2+q_2)(0)-(-p_1)(-p_1+p_2)=-p_1^2+p_1p_2,
\endalign
$$
so
$$
\{ \Cal R_1,\Cal R_2\} = -q_1(q_2) - q_2(-p_1^2+p_1p_2) = q_2\Cal R_2.
$$
Thus, condition (L1) holds, but the stronger condition (L2) is 
{\it not} satisfied in this case.

As an indication of the importance of the Lagrangian subvariety $V_M$
we shall discuss an observation of Givental and Kim (\cite{Gi-Ki})
concerning the case $M=F_{1,2,\dots,n-1}(\C^n)$. Lagrangian manifolds
appear naturally in the theory of completely integrable Hamiltonian
systems, and it turns out for this $M$  that 
$V_M$ is such an example --- the
integrable system being the \ll one-dimensional Toda lattice\rrr,
or 1DTL for short.  

For the case $M=F_{1,2}(\C^3)$, the observation will 
follow immediately from the
calculation of $\widetilde{QH}^\ast(F_{1,2}(\C^3);\C)$
in Appendix 1, once we have given the definition of the
1DTL.  To end this section, therefore, we shall sketch the
theory of the 1DTL for $n=3$. A more detailed but very elementary
discussion may be found in \cite{Gu2}.

The 1DTL is a system of first order o.d.e. in $2(n-1)$
functions of $t\in\R$.  For $n=3$ this system is
$$
\align
\dot a_1 &= a_1(b_1-b_2)\\
\dot a_2 &= a_2(b_2-b_3)\\
&\\
\dot b_1 &= -a_1\\
\dot b_2 &= -a_2+a_1\\
\dot b_3 &=\ \ \ \ \ \ \ \ \ a_2
\endalign
$$
where $b_1+b_2+b_3=0$ and $a_1,a_2> 0$.  With $u_i=\log a_i$
this can be written in the form
$$
\pmatrix
\ddot u_1 \\ \ddot u_2
\endpmatrix
=
\pmatrix 
-2 & 1\\
1 & -2
\endpmatrix
\pmatrix
e^{u_1} \\ e^{u_2}
\endpmatrix.
$$
However, a different matrix formulation  reveals
the interesting geometrical structure of the system, namely
the \ll Lax equation\rr
$$
\dot X = [X,Y]
$$
where
$$
X=
\pmatrix
b_1 & a_1 & 0\\
1 & b_2 & a_2\\
0 & 1 & b_3
\endpmatrix,
\quad
Y=
\pmatrix
0 & a_1 & 0\\
0 & 0 & a_2\\
0 & 0 & 0
\endpmatrix.
$$
If we define
$$
B=\{ (t_1,t_2) \in \R^2 \st t_1,t_2 > 0 \}
$$
we may identify $T^\ast B$ as the \ll phase space\rr
$$
T^\ast B \cong
\left\{
\pmatrix
s_1 & t_1 & 0\\
1 & s_2 & t_2\\
0 & 1 & s_3
\endpmatrix
\ \text{such that}\ 
(t_1,t_2)\in B, s_1,s_2,s_3\in\R, s_1+s_2+s_3=0\right\}
$$
of the system.  This is of course a (real) symplectic manifold.  It
can be shown that the o.d.e. $\dot X=[X,Y]$ describes the integral 
curves of a Hamiltonian vector field on $T^\ast B$.

Now, the form of the Lax equation implies that the solutions
$X$ are of the form $X(t)=A(t)V A(t)^{-1}$ for some function
$A:\R\to GL_3\R$.  Without loss of generality we have
$A(0)=I$ and so $X(0)=V$. In fact, it is possible to find $A$
explicitly\footnote{ Namely:  $A(t)$ is the matrix obtained by
applying the Gram-Schmidt orthogonalization process to the
columns of the matrix $\exp\,tV$. From this one obtains
expressions for $a_i(t)$ and $b_i(t)$ as rational functions
of exponential functions.}
and thus the solution of the original o.d.e., but
we shall not need the explicit formula here.  We just need the
observation that {\it the function $\det (X(t)+\la I)$ is independent
of $t$.} Indeed, we have $\det\, X(t)+\la I=\det \, V+\la I$ for any $\la$, hence
the coefficients of the powers of $\la$ are \ll conserved quantities\rr
along each solution $X(t)$.  These coefficients are as follows:
$$
\align
\la^3:\quad &1\\
\la^2:\quad &b_1+b_2+b_3 \ (=0)\\
\la^1:\quad &b_1b_2 +b_2b_3 + b_3b_1 - a_1 - a_2\  (= g(a,b), \ \text{say})\\
\la^0:\quad &b_1b_2b_3 -b_3a_1 -b_1a_2 \  (= h(a,b), \ \text{say})
\endalign
$$

In classical language, the two nontrivial conserved quantities
$g$ and $h$ are \ll first integrals\rr of the system, and they lead
to (another) method of finding the solution:  substitute $a_1=-\dot b_1$,
$a_2=\dot b_3$ into the equations
$$
\align
g(a,b)&=C_1\\
h(a,b)&=C_2
\endalign
$$
(for constants $C_1,C_2$)
and then try to solve the reduced system of first order
o.d.e for $b_1,b_2,b_3$.

In modern language, the equations $g(a,b)=C_1$ and
$h(a,b)=C_2$ define a \ll Lagrangian leaf\rr of
the foliation of $T^\ast B$ given by the corresponding 
Hamitonian vector fields
$H_g,H_h$. Each solution curve $X$ lies entirely within
such a leaf.

{\it Comparison with our earlier discussion of $M=F_{1,2}(\C^3)$
reveals that this Lagrangian leaf --- for $C_1=C_2=0$, and
in the complex variable case --- is precisely the manifold $V_M$.}

A complete proof of this surprising coincidence
in the case of $n\times n$ matrices --- which was sketched only
rather briefly in \cite{Gi-Ki} --- can be found in \cite{Ci} and \cite{Ki}.
The main point is to prove that the conserved quantities
of the 1DTL are relations in the quantum cohomology
ring of $F_{1,2,\dots,n-1}(\C^n)$.  The proofs given in these references all
use indirect methods (with a view to future generalizations).  
An elementary direct proof can be found in  \cite{Gu-Ot}, together with the
sketch of a generalization to the case of the (infinite dimensional)
\ll periodic flag manifold\rr and the periodic 1DTL.

\head
\S 6 A system of differential operators
\endhead

References:  \cite{Gi-Ki}, \cite{Au3}

We assume as in \S 5 that  $H^2(M;\Z)=\oplus_{i=1}^r \Z b_i$ 
generates $H^\ast(M;\Z)$.
Writing $t=\sum_{i=1}^r t_i b_i$,
we have in this situation  a generating function for the
(finitely many) products 
$b_1^{j_1}\dots b_r^{j_r}$ in  $H^\ast(M;\Z)$:
$$
v(t)=e^{t} = \sum_{l\ge 0} \frac 1{l!} t^l 
= \sum_{j,l \text{ s.t. } \sum j_i = l}\ 
\frac{t_1^{j_1}\dots t_r^{j_r}}{j_1!\dots j_r!} b_1^{j_1}\dots b_r^{j_r}.
$$
Evidently
$$
\left.\frac{\bd^{j_1+\dots+j_r}}
{\bd t_1^{j_1} \dots \bd t_r^{j_r}}\  v(t) \right |_{0}
=  b_1^{j_1}\dots b_r^{j_r},
$$
i.e. the products are given explicitly by the derivatives of the 
generating function $v$.
There is also a \ll scalar\rr generating function
$$
V(t)= \lan v(t), M \ran,
$$
which carries the same information as $v$.

\proclaim{Example: $M=\C P^n$}
\endproclaim

\no We have  $v(t)= 1 + tb + \frac 1{2!} t^2 b^2 +
\dots + \frac 1{n!} t^n b^n$, where $b=x_1$.
The scalar generating function is just $V(t)=\frac 1{n!} t^n$.

\proclaim{Example: $M=F_{1,2}(\C^3)$}
\endproclaim

\no Using the notation of Appendix 1, we have:
$v(t)= 1 + t_1 a + t_2 b +
\frac 1{2!}(t_1^2 a^2 + 2t_1t_2 ab + t_2^2 b^2)
+ \frac 1{3!}(3t_1^2 t_2 a^2b + 3t_1 t_2^2 ab^2)$.
The scalar generating function is  
$V(t)=\frac 1{3!}(3t_1^2 t_2 + 3t_1 t_2^2)$.

For quantum products we can make a similar definition.  Namely,
we define
$$
v(t,q)= \sum_{l\ge 0} \frac 1{l!}t\circ \dots\circ t
$$
where the indicated product $t\circ \dots\circ t$ is
the quantum product of $l$ copies of $t$. We also have
$$
V(t,q)= \lan v(t,q), M \ran 
$$
We shall not discuss the
convergence of this formal series.  Nevertheless, as in the case of
ordinary cohomology, we have (formally, at least), the following
equation for the quantum products:
$$
\left.\frac{\bd^{j_1+\dots+j_r}}{\bd t_1^{j_1} \dots \bd t_r^{j_r}}\  
v(t) \right |_{0}
=  (b_1\circ\dots \circ b_1)\circ\dots\circ (b_r\circ\dots\circ b_r),
$$
where $b_1,\dots,b_r$ appear (respectively)
$j_1,\dots,j_r$ times on the right hand side.

\proclaim{Example: $M=\C P^n$}
\endproclaim

\no We have
$v(t,q)=\sum_{l\ge 0} \frac 1{l!}t^l b\circ \dots\circ b$ (with $b=x_1$),
and $V(t,q)=\sum_{s\ge 0} \frac{t^{(n+1)s+n}}{((n+1)s+n)!}q^s$.

Now we introduce the differential operators referred to in the
title of this section.  

\proclaim{Definition}  For any $R\in S\otimes \La$,  we define
a differential operator $R^\ast$ by
\newline
$R^\ast = R(\frac{\bd}{\bd t_1},\dots,\frac{\bd}{\bd t_r},
q_1,q_1^{-1},\dots,q_r,q_r^{-1} )$.
\endproclaim

\no In other words, to produce $R^\ast$ from $R$, we identify the
coordinate function $p_i$ with the cohomology class $b_i$, and then 
regard this as the vector field $\bd/\bd t_i$ on $H^2(M;\C)$.
The following observation of \cite{Gi-Ki} is also explained in
\cite{Au3} (Theorem 3.3.1).

\proclaim{Theorem}  $R^\ast V(t,q)=0 \iff R\in \Cal I$.
\endproclaim

\demo{Proof} We have
$$
\align
R^\ast V = 0 &\iff 
\left.\frac{\bd^{j_1+\dots+j_r}}{\bd t_1^{j_1} \dots \bd t_r^{j_r}}\  
R^\ast V \right |_{0} = 0 \text{ for all } j\\
&\iff (p_1^{j_1}\dots p_r^{j_r} R)^\ast V \vert_0 
= 0 \text{ for all } j\\
&\iff \lan (p_1^{j_1}\dots p_r^{j_r} R)^\ast v, M \ran \vert_0 
= 0 \text{ for all } j\\
&\iff \lan b_1^{j_1}\dots b_r^{j_r} R, M \ran = 0 \text{ for all } j
\endalign
$$
where $b_1^{j_1}\dots b_r^{j_r}R$ is interpreted as an expression in
quantum cohomology, i.e. all products are quantum products.
By the nondegeneracy of the intersection form, the last
condition is equivalent to the vanishing of $R$ (considered
as a quantum polynomial), i.e. to the statement
$R \in \Cal I$.
\enddemo

This result is somewhat tautological, since the generating function
$V$ clearly contains by definition all the quantum products. However,
it gives an interesting new point of view on the quantum cohomology algebra:  
it is the algebra of differential operators which annihilate the 
generating function $V$.

\head
\S 7 The role of the flat connection
\endhead

References:  \cite{Gi-Ki}, \cite{Gi1}-\cite{Gi6}

In sections \S 4-\S 6 we discussed three \ll analytic\rr 
aspects of quantum cohomology.
The key to understanding these is the family of flat connections  
$\nabla\sla=d+\la \om$ of \S 1, with $\om_t(X)(Y)=X\circ_t Y$.

\no{\it Flatness of $\nabla\sla$ via a generating function.}
In \S 2 we stated that  $\nabla\sla$ is  a flat connection,
for each $\la$.  Flatness of $\nabla\sla$ (for all $\la$) is equivalent to
the conditions
$$
d \om = 0,\quad \om \wedge \om = 0
$$
(see \S 1).  The second condition follows from the commutativity
and associativity of $\circ_t$.  To establish the first condition,
it suffices to find a function $K:H^2(M;\C)\to \End(W)$ such that $\om = dK$.

Before proceeding with this, we remark that there is an obvious
generalization of the triple products $\ABC\subD$, namely
$$
\lan A_1\vert \dots \vert A_i \ran\subD\ =
\vert
\HolD^{\tilde A_1,p_1} \cap \dots \cap \HolD^{\tilde A_i,p_i} 
\vert
$$
for any $i\ge 3$.  It can be shown (with the appropriate technical
assumptions, cf. \S 3) that
$$
a_1\circ \dots \circ a_i = 
\sum_{D\in H_2(M;\Z)}(a_1\circ \dots \circ a_i)\subD q^D
$$
where
$$
\lan (a_1\circ \dots \circ a_i)\subD, B \ran =
\lan A_1\vert \dots \vert A_i\vert B \ran\subD
$$
for all homology classes $B$; hence these $i$-fold products are determined
\footnote{There is another kind of $i$-fold product 
in which the basepoints $p_1,\dots,p_i$
are allowed to vary (and which are of more
interest, from the point of view of enumerative geometry). 
These $i$-fold products are not 
in general determined by the triple products, and the situation is
more complicated. For $i=3$, both definitions coincide.} 
by the quantum products, and therefore by the $3$-fold (triple) products, for
$i\ge 3$. Using this notation, the generating function $V(t,q)$
of \S 6 is (ignoring the \ll quadratic part\rrr)
$$
V(t,q)= \sum_{l\ge 3; D} \frac 1 {l!} \lan T \vert \dots \vert T \ran\subD q^D
$$
where $T$ (the Poincar\acuteaccent e dual of $t=\sum_{i=1}^r t_i b_i$)
appears $l$ times in the indicated term.

A slight modification gives a definition of $2$-fold products, namely
$$
\lan A_1\vert A_2 \ran\subD\ =
\vert
( \HolD^{\tilde A_1,p_1} \cap \HolD^{\tilde A_2,p_2} )/\C^\ast
\vert
$$
where $\C^\ast$ is identified with the subgroup of linear
fractional transformations (of $\C P^1$) which fix $p_1$ and $p_2$.
These products are also related to the triple products, because of
the \ll divisor property\rr
$$
\lan A_1\vert A_2 \vert X \ran\subD = \lan A_1\vert A_2\ran\subD\, \lan x,D\ran
$$
where $x$ is any element of $H^2(M;\Z)$ and $D\ne 0$.  (Roughly
speaking, this formula follows from the fact that a nonconstant holomorphic map
$f:\C P^1\to M$ in the homotopy class $D$ must hit a subvariety $X$ of
complex codimension $1$ in $\lan x,D\ran$ points --- this is a topological 
necessity, beyond the influence of $A_1,A_2$.)

Let us return now to the function $K$.  We claim that the following explicit
formula defines a function with the desired property $\om = dK$:
$$
\lan K( t )(a), C\ran = \lan A\vert C\vert T\ran_0 +
\sum_{D\ne 0} \lan A\vert C\ran\subD e^{\lan t,D\ran}
$$
where $t\in H^2(M;\Z)$ and $a\in H^\ast(M;\Z)$
are fixed, and where $C$ is a general element of $H_\ast(M;\Z)$.
The function $K$ may be extended by complex linearity to the
case of  $t\in H^2(M;\C)$ and $a\in H^\ast(M;\C)$.

Having produced a candidate for $K$, 
the proof of the claim will be an elementary verification. We choose
$x\in H^2(M;\C)$ and regard it as a vector field on $H^2(M;\C)$. 
Using the differentiation formula $x e^{\lan t,D\ran}=
\lan x,D\ran e^{\lan t,D\ran}$, we obtain:
$$
\align
x\cdot\lan K(t)(a),C\ran &= 
x\cdot\lan act,M\ran +
\sum_{D\ne 0} \lan A\vert C\ran\subD\, x\cdot e^{\lan t,D\ran}\\
&=\lan acx,M\ran +
\sum_{D\ne 0} \lan A\vert C\ran\subD \lan x,D\ran e^{\lan t,D\ran}\\
&=\lan A\vert C\vert X\ran_0 + 
\sum_{D\ne 0} \lan A\vert C\vert X\ran\subD\,  e^{\lan t,D\ran}\\
&= \lan x\circ_t a,C\ran.  
\endalign
$$
Hence $dK(t)(a)(x) = x K(t)(a) =  x\circ_t a = \om_t(a)(x)$,
i.e. $dK = \om$, as required.

\no{\it Generating functions as solutions of integrable systems.}
The flatness of the connection $\nabla^{\la}=d+\la \om$ (on
the simply connected manifold $W$) implies that it is gauge-equivalent
via some gauge transformation $H$ to the trivial connection
$d$.  In other words, the equation $d\om + \la \om\wedge \om = 0$ implies that
$\la\om = H^{-1} dH$ for some map $H:W\to Gl(W)$.
(This is a standard fact from differential geometry.  Note that
$H$ will depend on $\la$ here.)
It can be shown that \ll generating function\rr $H$ characterizes the
quantum product in a similar way to the generating function
$V$ of \S 6. This fact was used (for example) in \cite{Ki} to establish a
relation between the quantum cohomology of the flag manifold
$G/B$ and the Toda lattice associated to a Lie group $G$, generalizing
the case $G=U_n$ which was discussed in \S 5. A brief
explanation of this argument can be found in \cite{Gi6}.

\no{\it Lagrangian subvarieties.}
Another application of the flat connection $\nabla^{\la}$ is the
proof that the variety $V_M$ of \S 5 is Lagrangian. We shall just
sketch the argument (of \cite{Gi-Ki}, \cite{Au4}) briefly.

The first ingredient is the subvariety $\Si_M$ of $T^\ast B$
defined by the \ll characteristic polynomials of the linear
transformations $\om_t(x) \in \End(W)$, for $x\in H^2(M;\C)$\rrr.
That is,
$$
\Si_M=\{ ([t],\la) \in B\times H^2(M;\C)^\ast \st
\det \om_t(x) - \la(x) I = 0
\text{ for all } x\in H^2(M;\C) \}.
$$
Let us assume that the linear transformation $\om_t(x)$
has $s+1 (=\dim W)$ distinct eigenvalues for at least
one (and hence almost all) values of $t$. Then it can be
shown that $\Si_M$ and $V_M$ {\it coincide.} The variety
$V_M$ therefore has $s+1$ branches, the $i$-th branch being given
by the equations
$$
p_1=\la_i(b_1),\ \ \dots,\ p_r=\la_i(b_r)
$$
where $\la_i$ is the $i$-th eigenvalue. (Note that the entries
of the matrix of $\om_t(b_i)$ with respect to the basis
$b_0,\dots,b_s$, and hence the eigenvalues
$\la_i(b_j)$,  are functions of $q_1,\dots,q_r$.)

Next, from the fact that $d\om=0$, it can be deduced that $d\la_i=0$.
The eigenvalue $\la_i$ can be regarded as a $1$-form on $B$,
and the $i$-th branch of $\Si_M$ is simply the image of this
$1$-form. But it is well known (see \cite{We}) that the
image of a closed $1$-form is Lagrangian, so this completes
the proof.

\no{\it Landau-Ginzburg potentials.}
So far we have not mentioned the 
Landau-Ginzburg potentials of \S 4. It turns out
that these are related to the mirror symmetry phenomenon.
As explained in \cite{Au2}, mirror symmetry can be viewed 
as a correspondence betweeen \ll Landau-Ginzburg models\rr 
(which produce Frobenius manifolds from singularity theory)
and \ll field theory models\rr (with Frobenius manifolds given
by quantum cohomology). We have seen that the quantum cohomology of
the Grassmannian is related to the singularities of the 
Landau-Ginzburg potential $\tilde P$. A mirror principle for
the flag manifold $F_n$ is proposed in \cite{Gi6}.

\head
Appendix 1: $QH^\ast(F_{1,2}(\C^3))$
\endhead

{\eightpoint
We shall calculate --- in obsessive detail --- the quantum 
cohomology of the flag manifold
$$
F_3=F_{1,2}(\C^3) =
\{ (L,V) \in Gr_1(\C^3) \times Gr_2(\C^3) \st L\sub V \}.
$$
This is a complex manifold of real dimension $6$ (so $n=3$
in the terminology of \S 2 and \S 3).

The following \ll Borel description\rr of $H^\ast(F_3;\Z)$ is
well known:
$$
H^\ast(F_3;\Z) \cong \Z[x_1,x_2,x_3]/\lan \si_1,\si_2,\si_3\ran
$$
where $\si_1,\si_2,\si_3$ are the elementary symmetric functions
of $x_1,x_2,x_3$.  Geometrically, $x_i=-c_1 \Cal L_i$, where
$\Cal L_1,\Cal L_2,\Cal L_3$ are the complex line bundles on $F_3$
whose fibres over $(L,V)$ are $L$, $L^\perp\cap V$, $V^\perp$
respectively.

We may choose additive generators as follows:
$$
\matrix
H^0(F_3;\Z) & H^2(F_3;\Z) & H^4(F_3;\Z) & H^6(F_3;\Z) \\
& & & \\
1 & x_1 & x_1^2 & x_1^2 x_2 \\
 & x_2 & x_1 x_2 & \\
\endmatrix
$$
The remaining cup products are determined by the relations
$x_1^2 + x_2^2 + x_1x_2 = 0$, 
$x_1x_2^2 + x_1^2x_2 = 0$
(which imply $x_1^3=0, x_2^3=0$).

This description ignores the complex manifold structure of $F_3$.  Since
we shall be considering holomorphic maps, it is more appropriate to
use the \ll Schubert description\rr of $H^\ast(F_3;\Z)$, which amounts
to replacing $x_1$, $x_2$ by $a=x_1$, $b=x_1+x_2$.  Geometrically,
$a=c_1 \Cal L^\ast$ and $b=c_1\Cal V^\ast$, where $\Cal L$, $\Cal V$
are the holomorphic \ll tautological\rr bundles whose fibres over
$(L,V)$ are $L$, $V$ respectively.  We then choose the following 
additive generators:
$$
\matrix
H^0(F_3;\Z) & H^2(F_3;\Z) & H^4(F_3;\Z) & H^6(F_3;\Z) \\
& & & \\
1 & a & a^2 & a^2b=ab^2 \\
 & b & b^2 & \\
\endmatrix
$$
This time the remaining cup products are determined by 
$ab=a^2+b^2$
(which implies $a^3=0, b^3=0$).
From the definition of $a$, $b$ we have $a^2=x_1^2$,
$b^2=x_1x_2$, and $a^2b=x_1^2x_2$, so the only difference
between the last table and the previous table is that $x_2$
has been replaced by $b$.  From the  table we have
$a^c=b^2$, $b^c=a^2$, $(a^2)^c=b$, $(b^2)^c=a$.

With respect to a fixed reference flag $E_1\sub E_2\sub \C^3$ there
are six Schubert varieties.  Excluding the trivial cases
\ll$F_3$\rr and \ll a point\rrr, we list these below, together
with their Poincar\acuteaccent e dual cohomology classes.

\no (1)\quad $\{L\sub E_2\} \ =\ \PD(a) \ \in \ H_4(F_3;\Z)$

\no (2)\quad  $\{E_1\sub V\} \ =\ \PD(b) \ \in \ H_4(F_3;\Z)$

\no (3)\quad  $\{L= E_1\} \ =\ \PD(a^2) \ \in \ H_2(F_3;\Z)$

\no (4)\quad  $\{V= E_2\} \ =\ \PD(b^2) \ \in \ H_2(F_3;\Z)$.

\no Of course, $\{L\sub E_2\}$ is an abbreviation for
$\{ (L,V)\in F_3 \st L\sub E_2 \}$, etc.

We are now ready to consider holomorphic maps
$f:\C P^1\to F_3$.  We begin by choosing the following basis for
$H_2(F_3;\Z)$:
$$
H_2(F_3;\Z)\ =\ \Z A^c \oplus \Z B^c\quad\quad
\text{(where $A^c=\PD(a^c)=\PD(b^2), B^c=\PD(b^c)=\PD(a^2)$).}
$$
If $[f]=d_1 A^c + d_2 B^c$ for some $d_1,d_2\in \Z$,
then it follows
\footnote{
For example, $f^\ast a=\lan a,d_1A^c+d_2 B^c\ran =
d_1\lan a,A^c\ran + d_2\lan a,B^c\ran = d_1$.
} that
$$
d_1=f^\ast a,\quad d_2 = f^\ast b
$$
(with respect to a fixed choice of generator of $H^2(\C P^1;\Z)$).
It can be shown that the space $\Hol^{F_3,p}_{d_1,d_2}$ is
nonempty if and only if either (a) $d_2\ge d_1\ge 0$ or
(b) $d_2=0, d_1\ge 0$.

With respect to the above choice of basis we can write
$q^{(d_1,d_2)}=q_1^{d_1}q_2^{d_2}$, where $q_1=q^{A^c}$
and $q_2=q^{B^c}$.  We have $\vert q_1\vert = \vert q_2\vert = 4$.
To calculate the quantum product $\circ$,
we begin by investigating the numerical condition
for $\lan X\vert Y\vert Z\ran\subD\ne 0$.  It is known that
$\lan c_1 TF_3, D\ran = 2d_1 + 2d_2$, so the numerical condition is
$$
\vert x\vert + \vert y\vert + \vert z\vert = 6 + 4d_1 + 4d_2.
$$
Since $\Hol^{F_3,p}_{d_1,d_2}$ is empty when either of $d_1$ or
$d_2$ is negative, we have
$$
x\circ y = \sum_{d_1,d_2\ge 0} (x\circ y)_{d_1,d_2}\ 
q_1^{d_1} q_2^{d_2}.
$$
Since the degree $\vert (x\circ y)_{d_1,d_2}\vert$ is given
by $\vert x\vert + \vert y\vert - 4d_1 -4d_2$, and this
must be $0$, $2$, $4$, or $6$, the relevant values of $(d_1,d_2)$
are severely restricted.
We shall calculate all possible quantum products of the additive basis 
elements $a$, $b$, $a^2$, $b^2$, $a^2b$.

\proclaim{Proposition 1} $a\circ a=a^2 + q_1$.
\endproclaim

\demo{Proof} We have $a\circ a= (a\circ a)_{0,0} + (a\circ a)_{1,0}q_1
+ (a\circ a)_{0,1}q_2$.  Now, $(a\circ a)_{0,0}$ is necessarily $a^2$,
so it remains to calculate the degree $0$ cohomology classes
$\la=(a\circ a)_{1,0}$ and $\mu=(a\circ a)_{0,1}$.

By definition, $\la = \lan A\vert A\vert Z\ran_{1,0}$, where $Z$
denotes the generator of $H_0(F_3;\Z)$.  We must look for three suitable
representatives of the homology classes $A,A,Z$.  We try

\no(a)\quad $\{L\sub E'_2\}$\quad (representing $A$)

\no(b)\quad $\{L\sub E''_2\}$\quad (also representing $A$)

\no(c)\quad $(E_1,E_2)$\quad (a single point of $F_3$)

\no where $E_1,E_2,E'_2,E''_2$ are to be chosen (if possible) so
that there are only finitely many holomorphic maps of degree
$(1,0)$ touching (a), (b) and (c) at three distinct points.

Now, any holomorphic map of degree $(1,0)$ is of the form
$\P(H)\to F_3$, $L\mapsto (L,H)$, where $H$ is a fixed two-dimensional
subspace of $\C^3$.  Therefore, a precise formulation of the
problem is as follows:  we must choose  $E_1,E_2,E'_2,E''_2$
so that there are only finitely many configurations
$(L',L'',L''',H)$ --- where $L', L'', L'''$ are distinct lines
in $H$ --- such that

\no(a)\quad $L'\sub E'_2$

\no(b)\quad $L''\sub E''_2$

\no(c)\quad $L''' = E_1, \ H=E_2$

\no Let us choose the two-dimensional subspaces $E_2,E'_2,E''_2$
in general position, i.e. such that the intersection of any two of
them is a line, and the intersection of all three is the origin.
Let us choose $E_1$ to be any line in $E_2$.  Then there is a
{\it single} configuration $(L',L'',L''',H)$ satisfying (a), (b)
and (c), namely
$$
L'=H\cap E'_2,\ L''=H\cap E''_2,\ L''' = E_1,\ H=E_2.
$$
We conclude that $\la=1$.

Next we calculate $\mu=\lan A\vert A\vert Z\ran_{0,1}$ by a similar
method.  The new feature is that holomorphic maps of degree $(0,1)$
are of the form $\P(\C^3/K)\to F_3$, $V/K\mapsto (K,V)$, where
$K$ is a fixed line in $\C^3$.  So we must choose $E_1,E_2,E'_2,E''_2$
so that there are only finitely many configurations
$(K,V',V'',V''')$ --- where $V',V'',V'''$ are distinct 
two-dimensional subspaces containing the line $K$ --- such that

\no(a)\quad $K\sub E'_2$

\no(b)\quad $K\sub E''_2$

\no(c)\quad $K = E_1, \ V'''=E_2$.

\no By choosing $E_2,E'_2,E''_2$ in general position, and $E_1\ne E'_2\cap E''_2$,
we see that there are {\it no} such configurations.  We conclude that
$\mu=0$.
\enddemo

\proclaim{Proposition 2} $b\circ b=b^2 + q_2$.
\endproclaim

\demo{Proof} Similar to the proof of Proposition 1.
\enddemo

\proclaim{Proposition 3} $a\circ b=ab$.
\endproclaim

\demo{Proof} We have to show that $(a\circ b)_{1,0}=(a\circ b)_{0,1}=0$.
First we have $(a\circ b)_{1,0}=\lan A\vert B\vert Z\ran_{1,0}$.  As
representatives of $A,B,Z$ we try

\no(a)\quad $\{L\sub E_2\}$\quad (representing $A$)

\no(b)\quad $\{E'_1\sub V\}$\quad (representing $B$)

\no(c)\quad $(E''_1,E''_2)$\quad (a single point of $F_3$)

\no We must choose $E'_1,E''_1,E_2,E''_2$ so that there are only
finitely many configurations
$(L',L'',L''',H)$ --- where $H$ is two-dimensional and
$L',L'',L'''$ are distinct lines in $H$ --- such that

\no(a)\quad $L'\sub E_2$

\no(b)\quad $E'_1\sub H$

\no(c)\quad $L''' = E''_1, \ H=E''_2$

\no Let us choose $E''_1,E_2,E''_2$ arbitrarily, and $E'_1$ such that
$E'_1\not\sub E''_2$.  Since (b) and (c) imply that $E'_1\sub E''_2$,
there are no such configurations.  Hence $(a\circ b)_{1,0}=0$.
A similar calculation shows that $(a\circ b)_{0,1}=0$.
\enddemo

\proclaim{Proposition 4} $a\circ b^2=ab^2$.
\endproclaim

\demo{Proof} We have $a\circ b^2= ab^2 + (a\circ b^2)_{1,0}q_1
+ (a\circ b^2)_{0,1}q_2$.

Let $(a\circ b^2)_{1,0}=\la a+\mu b$.
We have $\la = \lan A\vert B^2\vert B^2\ran_{1,0}$. As
representatives of the homology classes $A,B^2,B^2$, we try:

\no(a)\quad $\{L\sub E_2\}$

\no(b)\quad $\{V=E'_2\}$

\no(c)\quad $\{V=E''_2\}$

\no We will choose $E_2,E'_2,E''_2$ so
that there are only finitely many configurations
$(L',L'',L''',H)$ --- where $L', L'', L'''$ are distinct lines
in $H$ --- such that

\no(a)\quad $L'\sub E_2$

\no(b)\quad $H= E'_2$

\no(c)\quad $H=E''_2$

\no Thus $E'_2=E''_2$, but this is impossible if we choose
$E_2,E'_2,E''_2$ in general position.  Hence $\la=0$.

Next, we have $\mu = \lan A\vert B^2\vert A^2\ran_{1,0}$.  
As representatives of the homology classes $A,B^2,A^2$, we try:

\no(a)\quad $\{L\sub E_2\}$

\no(b)\quad $\{V=E'_2\}$

\no(c)\quad $\{L=E''_1\}$

\no We will choose $E_2,E'_2,E''_1$ so
that there are only finitely many configurations
$(L',L'',L''',H)$ --- where $L', L'', L'''$ are distinct lines
in $H$ --- such that

\no(a)\quad $L'\sub E_2$

\no(b)\quad $H= E'_2$

\no(c)\quad $L'''=E''_1$

\no Thus, $E''_1\sub E'_2$, but we may choose $E_2,E'_2,E''_2$
so that this is false.  Hence $\mu=0$.

Now we turn to $(a\circ b^2)_{0,1}=\la a+\mu b$.
We have $\la = \lan A\vert B^2\vert B^2\ran_{0,1}$. As
representatives of the homology classes $A,B^2,B^2$, we try:

\no(a)\quad $\{L\sub E_2\}$

\no(b)\quad $\{V=E'_2\}$

\no(c)\quad $\{V=E''_2\}$

\no We will choose $E_2,E'_2,E''_2$ so
that there are only finitely many configurations
$(K,V',V'',V''')$ --- where $V',V'',V'''$ are distinct two-dimensional
subspaces containing the line $K$ --- such that

\no(a)\quad $K\sub E_2$

\no(b)\quad $V''= E'_2$

\no(c)\quad $V'''=E''_2$

\no Thus $K\sub E_2,E'_2,E''_2$, but this is impossible if we choose
$E_2,E'_2,E''_2$ in general position.  Hence $\la=0$.

Next, we have $\mu = \lan A\vert B^2\vert A^2\ran_{0,1}$.  
As representatives of the homology classes $A,B^2,A^2$, we try:

\no(a)\quad $\{L\sub E_2\}$

\no(b)\quad $\{V=E'_2\}$

\no(c)\quad $\{L=E''_1\}$

\no We will choose $E_2,E'_2,E''_1$ so
that there are only finitely many configurations
$(K,V',V'',V''')$ --- where $V',V'',V'''$ are distinct two-dimensional
subspaces containing the line $K$ --- such that

\no(a)\quad $K\sub E_2$

\no(b)\quad $V''= E'_2$

\no(c)\quad $K=E''_1$

\no Thus, $E''_1\sub E_2$, but we may choose $E_2,E'_2,E''_2$
so that this is false.  Hence $\mu=0$.
\enddemo

\proclaim{Proposition 5} $a^2\circ b=a^2b$.
\endproclaim

\demo{Proof} Similar to the proof of Proposition 4.
\enddemo

\proclaim{Proposition 6} $a\circ a^2=bq_1$.
\endproclaim

\demo{Proof} We have $a\circ a^2=  (a\circ a^2)_{1,0}q_1
+ (a\circ a^2)_{0,1}q_2$.

Let $(a\circ a^2)_{1,0}=\la a+\mu b$.
We have $\la = \lan A\vert A^2\vert B^2\ran_{1,0}$,
which is zero by the proof of Proposition 4.

Next, we have $\mu = \lan A\vert A^2\vert A^2\ran_{1,0}$.  
As representatives of the homology classes $A,A^2,A^2$, we try:

\no(a)\quad $\{L\sub E_2\}$

\no(b)\quad $\{L=E'_1\}$

\no(c)\quad $\{L=E''_1\}$

\no We will choose $E_2,E'_1,E''_1$ so
that there are only finitely many configurations
$(L',L'',L''',H)$ --- where $L', L'', L'''$ are distinct lines
in $H$ --- such that

\no(a)\quad $L'\sub E_2$

\no(b)\quad $L''= E'_1$

\no(c)\quad $L'''=E''_1$

\no Thus, $H=E'_1\oplus E''_1$ (if we choose $E'_1\ne E''_1$),
and $L'=E_2\cap  E'_1\oplus E''_1$ (if we 
choose $E_2\ne E'_1\oplus E''_1$), so the
configuration is determined uniquely.  Hence $\mu=1$.

Now we turn to $(a\circ a^2)_{0,1}=\la a+\mu b$.
We have $\la = \lan A\vert A^2\vert B^2\ran_{0,1}$, which is
zero by the proof of Proposition 4.

Next, we have $\mu = \lan A\vert A^2\vert A^2\ran_{0,1}$.  
As representatives of the homology classes $A,A^2,A^2$, we try:

\no(a)\quad $\{L\sub E_2\}$

\no(b)\quad $\{L=E'_1\}$

\no(c)\quad $\{L=E''_1\}$

\no We will choose $E_2,E'_1,E''_1$ so
that there are only finitely many configurations
$(K,V',V'',V''')$ --- where $V',V'',V'''$ are distinct two-dimensional
subspaces containing the line $K$ --- such that

\no(a)\quad $K\sub E_2$

\no(b)\quad $K= E'_1$

\no(c)\quad $K=E''_1$

\no Thus, $E'_1= E''_1$, but we may choose $E_2,E'_1,E''_1$
so that this is false.  Hence $\mu=0$.
\enddemo

\proclaim{Proposition 7} $b\circ b^2=aq_2$.
\endproclaim

\demo{Proof} Similar to the proof of Proposition 6.
\enddemo

This completes the computation of all quantum products of
degree at most $6$.
We shall be able to evaluate the other quantum products without doing any
further calculations of triple products.

\proclaim{Corollary 8} (i) $a\circ a^2b=b^2q_1 + q_1q_2$.
(ii) $b\circ a^2b=a^2q_2 + q_1q_2$.
\endproclaim

\demo{Proof} (i) $a\circ a^2b=a\circ(a^2\circ b) = 
(a\circ a^2)\circ b = (bq_1)\circ b
=b\circ b q_1 = b^2q_1 + q_1q_2$. (ii) Similarly.
\enddemo

\proclaim{Corollary 9} (i) $a^2\circ a^2=b^2q_1$.
(ii) $b^2\circ b^2=a^2q_2$.
(iii) $a^2\circ b^2=q_1q_2$.
\endproclaim

\demo{Proof} (i) $a^2\circ a^2=(a\circ a - q_1)\circ a^2 =
a\circ a\circ a^2 - a^2 q_1 = a\circ (bq_1) - a^2q_1 =
abq_1 - a^2q_1 = b^2q_1.$ (ii), (iii) Similarly.
\enddemo

\proclaim{Corollary 10} (i) $a^2\circ a^2b=aq_1q_2$.
(ii) $b^2\circ a^2b=bq_1q_2$.
\endproclaim

\demo{Proof} (i) $a^2\circ a^2b = a^2\circ (a^2\circ b) = (a^2\circ a^2)\circ b
= b^2q_1\circ b = (b^2\circ b)q_1 = aq_1q_2$.
(ii) Similarly.
\enddemo

\proclaim{Corollary 11} $a^2b\circ a^2b=abq_1q_2$.
\endproclaim

\demo{Proof} $a^2b\circ a^2b= (a^2\circ b)\circ a^2b
= a^2\circ(a^2q_2 + q_1q_2) = b^2q_1q_2 + a^2q_1q_2 =
abq_1q_2$.
\enddemo

What are the relations defining the quantum cohomology algebra?
By Theorem 2.2 of \cite{Si-Ti},
these are the \ll quantum versions\rr
of the Borel relations
$$
x_1 x_2 + x_2 x_3 + x_3 x_1 = 0,\quad x_1 x_2 x_3 = 0.
$$

We obtain:

$$
\align
x_1\circ x_2 + x_2\circ x_3 + x_3\circ x_1 &= 
a\circ (b-a) + (b-a)\circ (-b) + (-b)\circ a\\
&= -a\circ a - b\circ b + a\circ b\\
&=-q_1-q_2
\endalign
$$
{}
$$
\align
x_1\circ x_2\circ x_3 &=
a\circ (b-a)\circ (-b)\\
&=-a\circ b\circ b + a\circ a\circ b\\
&=-ab\circ b + a\circ ab\\
&=-(a^2+b^2)\circ b + a\circ (a^2+b^2)\\
&= bq_1-aq_2
\endalign
$$

\no Restricting attention to nonnegative powers of $q_1,q_2$,
we therefore have
$$
\widetilde{QH}^\ast(F_3;\C) \cong
\frac{\C[x_1,x_2,x_3,q_1,q_2]}
{\lan x_1x_2+x_2x_3+x_3x_1+q_1+q_2,
x_1x_2x_3+x_3q_1+x_1q_2 \ran}.
$$
It should be noted that this information --- the isomorphism
type of $\widetilde{QH}^\ast(F_3;\C)$ (or ${QH}^\ast(F_3;\C)$) ---
does not determine all the quantum products.
}

\head
Appendix 2: $QH^\ast(\Si_k)$
\endhead

{\eightpoint
We shall investigate the Hirzebruch surfaces 
$\Si_k = \P(\Cal O(0)\oplus \Cal O(-k) )$, where $\Cal O(i)$
denotes the holomorphic line bundle on $\C P^1$ 
with first Chern class $i$. Since
$$
\P(\Cal O(0)\oplus \Cal O(-k) ) \cong
\P(  \Cal O(k)\otimes (\Cal O(0)\oplus \Cal O(-k)) ) \cong
\P(\Cal O(k)\oplus \Cal O(0) )
$$
(as complex manifolds), we have $\Si_k\cong \Si_{-k}$.  Because
of this, we shall assume from now on that $k\ge 0$.
We shall use the following explicit description of $\Si_k$
(see \cite{Hi}):
$$
\Si_k =
\{ ([z_0;z_1;z_2],[w_1;w_2]) \in \C P^2 \times \C P^1 \st
z_1 w_1^k = z_2 w_2^k \}
$$
The following four subvarieties of $\Si_k$ will play an
important role in our calculations.
$$
\align
X_1 &= \{ z_2 = w_1 = 0 \} \\
X_2 &= \{ z_1 = z_2 = 0 \} \\
X_3 &= \{ z_1 = w_2 = 0 \} \\
X_4 &= \{ z_0 = 0 \} 
\endalign
$$
If $\Si_k$ is regarded as \ll$\Cal O(-k) \cup \text{ $\infty$-section}$\rrr, then
$X_1$ and $X_3$ are fibres, $X_2$ is the $0$-section, and
$X_4$ is the $\infty$-section.  We shall also denote by $X_i$ the
homology class in $H_2(\Si_k;\Z)$ represented by the
variety $X_i$.  As usual, the Poincar\acuteaccent e dual cohomology
class in $H^2(\Si_k;\Z)$ will be denoted $x_i$.

\proclaim{Proposition 1} The relations between the cohomology classes
$x_1,x_2,x_3,x_4$ are:

\no(1) $x_1=x_3$, $x_4=x_2+kx_1$.

\no(2) $x_1 x_3 = x_2 x_4 = 0$, $x_1 x_2 = x_1 x_4 = x_2 x_3 = x_3 x_4 = z$,
where $z$ is a generator of $H^4(\Si_k;\Z) \cong \Z$.

\no(3) $x_1^2 (=x_3^2) = 0$, $x_2^2=-kz$, $x_4^2 = kz$.
\endproclaim

\demo{Proof} (1)  follows from consideration of the meromorphic
functions on $\Si_k$ defined by the expressions $w_1/w_2$, $z_1/z_0$.
(2)  follows from consideration of intersections.  (3)  is a 
consequence of (1) and (2).
\enddemo

The cohomology algebra of $\Si_k$ is generated by
$x_1$ and $x_2$, by the Leray-Hirsch Theorem (Theorem 5.11 
and (20.7) of \cite{Bo-Tu}). From this, and by using the proposition, we obtain:
$$
\align
H^\ast(\Si_k;\Z) &\cong
\Z[x_1,x_2,x_3,x_4]/
\lan x_1-x_3,x_4-x_2-kx_1\ran +
\lan x_1 x_3, x_2 x_4 \ran \\
&\cong
\Z[x_1,x_4]/
\lan x_1^2, x_4^2-kz\ran.
\endalign
$$
The particular form of the first description --- where we
have written the ideal of relations as a sum of two
separate ideals --- is chosen to be consistent with the
usual description of the cohomology algebra of a toric
variety (see \cite{Od}).

It can be verified (from this description) that 
\footnote{
A suitable isomorphism is given by
$x_2\mapsto x_2+x_1$, $x_4\mapsto x_4-x_1$.}
$H^\ast(\Si_k;\Z) \cong H^\ast(\Si_{k+2};\Z)$, and that
$H^\ast(\Si_k;\Z) \not\cong H^\ast(\Si_{k+1};\Z)$.
This is consistent with the results of \cite{Hi},
where it was shown that $\Si_k$ is homeomorphic to
$\Si_l$ if and only if $k-l$ is even.

We may choose additive generators of $H^\ast(\Si_k;\Z)$ as follows:
$$
\matrix
H^0(\Si_k;\Z) & H^2(\Si_k;\Z) & H^4(\Si_k;\Z)\\
& & \\
1 & x_1 & z=x_1 x_4\\
 & x_4 &\\
\endmatrix
$$
From the proposition we have:
$$
x_1^c = x_2 = x_4 - kx_1, \ x_4^c = x_3 = x_1.
$$
We shall make use of the following representatives of
homology classes of $\Si_k$.  These are analogous to the Schubert
varieties in flag manifolds and Grassmannians, but the
situation is more rigid now:
$$
\align
X_1: &\text{\quad represented by}\ X_1,X_3 \ \text{or any \ll fibre\rr}\\
X_2: &\text{\quad represented only by}\ X_2\\
X_4: &\text{\quad represented only by}\ X_4
\endalign
$$

We turn now to holomorphic maps $f:\C P^1\to \Si_k$.  Such maps have
a very explicit description in terms of polynomials.  Namely, $f$ is of
the form
$$
f=([p_4;p_2 p_3^k; p_2 p_1^k],[p_1;p_3])
$$
where $p_1,p_2,p_3,p_4$ are arbitrary complex polynomials such that
$p_1,p_3$ have no common factor, and
$p_2,p_4$ have no common factor.  (The notation is chosen
so that $p_i\equiv 0$ if and only if $f(\C P^1) \sub X_i$.)

We define the \ll degree\rr of $f$ by $\deg f = (d,e)$,
where
$$
\align
d&= \max\{ \deg p_4, \deg p_2 p_3^k, \deg p_2 p_1^k \}\\
e&= \max\{ \deg p_1, \deg p_3 \}
\endalign
$$

\proclaim{Proposition 2}

\no(1) $\deg f = (d,e)$ if and only if $[f]= d X_1 + e X_2$
$(\in \pi_2 (\Si_k) \cong H_2(\Si_k;\Z) )$.

\no(2) There exists a holomorphic map $f$ such that 
 $\deg f = (d,e)$ if and only if either

(a) $d\ge ke \ge 0$, or 

(b) $d=0$, $e>0$.
\endproclaim

\demo{Proof} (1) Let $H_1,H_2$ be the restrictions to $\Si_k$ of the
tautologous line bundles on $\C P^1$, $\C P^2$.  Then the condition
$\deg f = (d,e)$ means that $d=-c_1 f^\ast H_2$ and
$e=-c_1 f^\ast H_1$.   It may be verified that
$c_1 H_1 = -x_1$ and $c_1 H_2 = -x_4$ (by using the
fact that the first Chern class of the dual of the tautologous
line bundle on $\C P^n$ is Poincar\acuteaccent e
dual to a hyperplane in $\C P^n$).  Hence the condition
$\deg f = (d,e)$ means that $d= f^\ast x_4$ and
$e= f^\ast x_1$. Let us write $[f]= u X_1 + v X_2$.  Then we have
$$
u = \lan x_4, uX_1 + v X_2 \ran 
= \lan x_4, f_\ast [\C P^1] \ran 
= \lan f^\ast x_4, [\C P^1] \ran 
=d
$$
and similarly $v=e$.  This completes the proof of (1).

\no (2)  This is obvious from the polynomial
expression for $f$.
\enddemo

We need one more ingredient in order to compute the
quantum cohomology of $\Si_k$:

\proclaim{Proposition 3}  $c_1 T\Si_k = 2x_4 - (k-2)x_1$.
\endproclaim

\demo{Proof} See section 3.3 of \cite{Od}, for example.
\enddemo

Let us choose $X_1$, $X_2$ as a basis of $H_2(\Si_k;\Z)$, and let
us denote the corresponding quantum parameters by $q_1,q_2$.  By
the usual definition we have
$$
\align 
\vert q_1 \vert &= 2\lan c_1 T\Si_k, X_1\ran = 
2\lan  2x_4 - (k-2)x_1, X_1\ran = 4\\
\vert q_2 \vert &= 2\lan c_1 T\Si_k, X_2\ran = 
2\lan  2x_4 - (k-2)x_1, X_2\ran = 2(2-k)
\endalign
$$

\no {\it Warning:}  The usual definition of $\vert q_1\vert$
and $\vert q_2\vert$ is designed to indicate the \ll expected\rr
dimension of $\Hol_{d,e}(\C P^1,\Si_k)$.  In case (a), i.e.
$d\ge ke \ge 0$, this is correct:
$$
\align
&\dim_{\C} \Hol_{d,e}(\C P^1,\Si_k) = 2d + (2-k)e + 2
\ \text{(from the polynomial description)}\\
&\dim_{\C} \Si_k + \lan c_1 T\Si_k, dX_1 + eX_2\ran = 2d + (2-k)e + 2.
\endalign
$$
In case (b), however, it is not correct in general.  Here we have
$d=0$ and $e>0$, so --- if $k>0$ --- the image of $f$ must be contained
in $X_2$. We obtain:
$$
\align
&\dim_{\C} \Hol_{d,e}(\C P^1,\Si_k) = 2e+1
\ \text{(from the polynomial description)}\\
&\dim_{\C} \Si_k + \lan c_1 T\Si_k, dX_1 + eX_2\ran = (2-k)e + 2.
\endalign
$$
This means that we cannot rely on the \ll numerical
condition\rr when performing our calculations. More precisely, the
problem is that $\Si_k$ is Fano only for $k=0,1$ and convex
only for $k=0$ --- see
the remarks at the end of \S 3 and at the end of this section.
Because of this, our calculations from now on will be merely heuristic.
However, we shall arrive at the correct answer, and some references
for a rigorous approach will be given later.

In calculating the triple products $\lan A\vert B\vert C\ran\subD$,
we shall only need the cases $D=(0,0)$, $(1,0)$, or $(0,1)$;
this follows from the polynomial representation of holomorphic
maps.  The polynomial representation shows also that

\no(1) the holomorphic maps of degree $(1,0)$ are precisely
the \ll fibres\rrr, and

\no(2) the holomorphic maps of degree $(0,1)$ are given by
the \ll$0$-section\rrr, $X_2$.

To calculate $\lan A\vert B\vert C\ran\subD$, we need to find representatives
of the homology classes $A,B,C$ such that there is only a finite
number of configurations $(P_1,P_2,P_3,\C P^1)$ such that
$P_1\in A$, $P_2\in B$, $P_3\in C$, where $\C P^1$ represents the
homotopy class $D$ and $P_1,P_2,P_3$ are three distinct points
in this $\C P^1$.  If such representatives do not exist, then
we shall say that $\lan A\vert B\vert C\ran\subD$ is \ll not computable\rrr.
If, for fixed $A$ and $B$, the necessary triple products 
$\lan A\vert B\vert C\ran\subD$ are
not computable, we shall say that 
the quantum product $a\circ b$ is \ll not directly computable\rrr.
(It may still be possible to compute $a\circ b$ indirectly, by expressing
it in terms of other known quantum products --- or, of course, by using
more general representatives of homology classes than we have so far allowed
ourselves.)

\no {\it Calculations for $\Si_0 = \C P^1 \times \C P^1$}

In this case, we can rely on the numerical condition.  We have
$\vert q_1\vert = \vert q_2\vert = 4$, and $x_1=x_3$, $x_2=x_4$.
We obtain
$$
\align
x_1\circ x_1 &= q_2\\
x_4\circ x_4 &= q_1\\
x_1\circ x_4 &= x_1x_4.
\endalign
$$
The quantum cohomology algebra is
$$
\widetilde{QH}^\ast(\Si_0;\C) \cong
\frac{\C[x_1,x_4,q_1,q_2]}
{\lan  x_1^2 - q_2,  x_4^2 - q_1  \ran}.
$$
This is just the tensor product of two copies of the
quantum cohomology algebra of $\C P^1$.

\no {\it Calculations for $\Si_1$}

In this case too, we can rely on the numerical condition.  The
problematical case is $D=(0,1)$, i.e. $d=0, e=1$, but in fact
the space of holomorphic maps of degree $D$ has the correct
dimension, as $2e+1 = e+2$ in this case.

We take $\vert q_1\vert = 4$ and $\vert q_2 \vert = 2$, and
proceed in the usual way. We shall try to calculate the 
six possible quantum products
involving the cohomology classes $x_1 (=x_3),x_2,x_4$.

\proclaim{Proposition 4} $x_1 \circ x_1 = x_2 q_2$.
\endproclaim

\demo{Proof} Let us write
$$
x_1 \circ x_1 = x_1^2 + \la q_1 + \mu q_2
$$
where $\la \in H^0(\Si_1;\Z)$ and $\mu \in H^2(\Si_1;\Z)$.

We have $\la = \lan X_1 \vert X_1 \vert Z\ran_{1,0}$, where
$Z$ is the generator of $H_0(\Si_1;\Z)$.  Let us choose
$X_1$, $X_3$, and any point of $\Si_1$ as representatives of
the three homology classes.  Then each holomorphic map of
degree $(1,0)$ (i.e. each \ll fibre\rrr) fails to intersect
all three representatives.  Hence $\la = 0$.

To calculate $\mu$, we must calculate 
$\lan \mu, Y\ran =\lan X_1 \vert X_1 \vert Y\ran_{0,1}$,
for two independent homology classes $Y\in H_2(\Si_1;\Z)$.

The triple product $\lan X_1 \vert X_1 \vert X_1\ran_{0,1}$
is equal to $1$, because we may choose three distinct 
\ll fibres\rr representing $X_1$, and then there is a unique
configuration $(P_1,P_2,P_3,X_2)$ which intersects these
fibres (respectively) in the points $P_1,P_2,P_3$.

The triple product $\lan X_1 \vert X_1 \vert X_2\ran_{0,1}$
is not computable, because there are infinitely many
configurations $(P_1,P_2,P_3,X_2)$ with the appropriate
intersection property, whatever representatives of
$X_1,X_1,X_2$ are chosen.

The triple product $\lan X_1 \vert X_1 \vert X_4\ran_{0,1}$
is equal to zero, because $X_2\cap X_4 = \es$.

We conclude that $\lan \mu, X_1\ran = 1$ and $\lan \mu, X_4\ran = 0$,
and hence $\mu = x_2$.
\enddemo

\proclaim{Proposition 5} $x_2 \circ x_2$ is not directly computable.
\endproclaim

\demo{Proof} (There is only one representative of the homology
class $X_2$.)
\enddemo

\proclaim{Proposition 6} $x_4 \circ x_4$ is not directly computable.
\endproclaim

\demo{Proof} (There is only one representative of the homology
class $X_4$.)
\enddemo

\proclaim{Proposition 7} $x_1 \circ x_4=x_1 x_4$.
\endproclaim

\demo{Proof} Let us write $x_1 \circ x_4=x_1 x_4 + \la q_1 + \mu q_2$.
We have $\la = \lan X_1 \vert X_4 \vert Z\ran_{1,0}$. By taking $Z$ as
any point in the complement of $X_1\cup X_4$, we see that $\la = 0$. 
A holomorphic map of degree $(0,1)$ cannot intersect
$X_4$, so $\mu=0$.
\enddemo   

\proclaim{Proposition 8} $x_1 \circ x_2$ is not directly computable.
\endproclaim

\demo{Proof} (The triple products $\lan X_1 \vert X_2 \vert Y \ran_{0,1}$ are
not computable, for any $Y$.)
\enddemo

\proclaim{Proposition 9} $x_2 \circ x_4=q_1$.
\endproclaim

\demo{Proof} Let us write $x_2 \circ x_4=x_2 x_4 + \la q_1 + \mu q_2$.
We have  $\la = \lan X_2 \vert X_4 \vert Z\ran_{1,0}$, and this is equal
to $1$, as there is a unique \ll fibre\rr which intersects $X_2$, $X_4$
and a point in the complement of $X_2\cup X_4$.
The triple products $\lan X_2 \vert X_4 \vert Y\ran_{0,1}$ are
necessarily zero, so we have $\mu=0$.
\enddemo

To summarize, we have now computed the following quantum products for
$\Si_1$:
$$
\align
x_1 \circ x_1 &= x_2 q_2\\
x_1 \circ x_4&=z\\
x_2 \circ x_4&=q_1.
\endalign
$$
From the relation $x_4=x_2+x_1$, we obtain the three quantum products
which were not directly computable:
$$
\align
x_4 \circ x_4 &= z + q_1\\
x_1 \circ x_2&=z - x_2q_2\\
x_2 \circ x_2&=-z + q_1 + x_2 q_2.
\endalign
$$
The quantum cohomology algebra is therefore
$$
\widetilde{QH}^\ast(\Si_1;\C) \cong
\frac{\C[x_1,x_4,q_1,q_2]}
{\lan  x_1^2 - x_2 q_2,  x_4^2 - z - q_1  \ran}.
$$

Hirzebruch surfaces are examples of toric manifolds, 
and it is possible to describe rational curves in general toric manifolds quite
explicitly (see \cite{Gu1}). A variant of quantum cohomology of toric manifolds
was introduced in \cite{Ba}, and its relation with the usual quantum cohomology
algebra was discussed in \cite{Co-Ka}  (examples 8.1.2.2 and 11.2.5.2)
and in \cite{Si} (section 4).
Unfortunately the situation is quite complicated because toric manifolds
are not convex in general (cf. the remarks at the end of \S 3). We have
already seen this in the case of the Hirzebruch surface $\Si_k$, which is
convex if and only if $k=0$, i.e. if and only if $\Si_k$ is
homogeneous. An alternative way of proceeding in the non-convex case
is given in section 15 of \cite{Au2}; the calculations for $\Si_1$
(or rather, for the isomorphic space $\Si_{-1}$) are given in section 16
of that article.  A rigorous treatment of the quantum cohomology of
$\Si_1$ --- and more generally, of the manifold obtained by
blowing up $\C P^n$ at a point --- appears in \cite{Ga}.
}

\eightpoint \it

\no  Department of Mathematics
\newline
Graduate School of Science
\newline
Tokyo Metropolitan University
\newline
Minami-Ohsawa 1-1, Hachioji-shi
\newline
Tokyo 192-0397, Japan

\no martin\@miyako.math.metro-u.ac.jp


\newpage
\Refs   
     
\widestnumber\key{Fu-Pa}

\ref 
\key  Au1
\by M. Audin
\paper Cohomologie quantique
\jour Ast\acuteaccent erisque
\vol 241
\paperinfo S\acuteaccent eminaire Bourbaki 806
\yr 1997
\pages 29--58
\endref

\ref 
\key  Au2
\by M. Audin
\paper An introduction to Frobenius manifolds, moduli spaces
of stable maps and quantum cohomology
\paperinfo http://irmasrv1.u-strasbg.fr/\~{}maudin/recherche.html
\jour
\yr 
\vol 
\pages 
\endref

\ref 
\key  Au3
\by M. Audin
\paper Symplectic geometry and volumes of moduli spaces
in quantum cohomology
\paperinfo preprint
\jour
\yr 
\vol 
\pages 
\endref

\ref 
\key  Au4
\by M. Audin
\paper Symplectic geometry in Frobenius manifolds
and quantum cohomology
\paperinfo (shortened version of \cite{Au3})
\jour J. Geom. and Physics
\yr 1998
\vol 25
\pages 183--204
\endref

\ref 
\key  Ba
\by V. Batyrev
\paper Quantum cohomology rings of toric manifolds
\jour Ast\acuteaccent erisque
(Journ\acuteaccent ees de g\acuteaccent eom\acuteaccent etrie 
alg\acuteaccent ebrique D'Orsay)
\yr 1993
\vol 218
\pages 9--34
\endref

\ref 
\key  Be1
\by A. Bertram
\paper Towards a Schubert calculus for maps from a Riemann surface
to a Grassmannian
\jour Int. Jour. of Math.
\yr 1994
\vol 5
\pages 811--825
\endref

\ref
\key  Be2
\by  A. Bertram
\paper Computing Schubert's calculus with Severi residues:
an introduction to quantum cohomology
\inbook   Moduli of Vector Bundles
\bookinfo Lec. Notes in Pure and App. Math. 179
\ed M. Maruyama
\publ Dekker
\yr 1996
\pages 1--10
\endref

\ref
\key  Be3
\by A. Bertram
\paper Quantum Schubert calculus
\jour Advances in Math.
\yr 1997
\vol 128
\pages 289--305
\endref

\ref
\key BCPP
\by G. Bini, C. de Concini, M. Polito and C. Procesi
\book On the work of Givental relative to mirror symmetry
\bookinfo Appunti dei Corsi Tenuti da Docenti della Scuola 
\publ Scuola Normale Superiore, Pisa
\yr 1998 (math.AG/9805097)
\endref

\ref
\key Bo-Tu
\by R. Bott and L.W. Tu
\book Differential Forms in Algebraic Topology
\publ Springer
\bookinfo Graduate Texts in Math. 82, 1982
\endref

\ref 
\key  Ci
\by I. Ciocan-Fontanine
\paper Quantum cohomology of flag varieties
\jour Int. Math. Res. Notes
\yr 1995
\vol 6
\pages 263--277
\endref

\ref
\key Co-Ka
\by D.A. Cox and S. Katz
\book Mirror Symmetry and Algebraic Geometry
\publ Amer. Math. Soc.
\yr 1999
\bookinfo Math. Surveys and Monographs 68
\endref

\ref
\key  Du
\by B. Dubrovin
\paper Integrable systems and 
classification of $2$-dimensional topological field theories
\inbook Integrable systems (Luminy, 1991)
\bookinfo Progr. Math. 115
\yr 1993
\publ  Birkh\"auser
\eds O. Babelon, P. Cartier and Y. Kosmann-Schwarzbach
\pages 313--359
\endref

\ref
\key Fu
\by W. Fulton
\book  Young Tableaux
\yr 1997
\publ Cambridge Univ. Press
\endref

\ref
\key Fu-Pa
\by W. Fulton and R. Pandharipande
\paper Notes on stable maps and quantum cohomology
\inbook Algebraic geometry---Santa Cruz 1995
\bookinfo Proc. Sympos. Pure Math., 62, Part 2
\eds J. Kollar, R. Lazarsfeld, and D.R. Morrison
\yr 1997
\pages 45--96
\publ Amer. Math. Soc.
\endref

\ref 
\key Ga
\by  A. Gathmann
\paper Counting rational curves with multiple points and
Gromov-Witten invariants of blow-ups
\jour 
\yr
\vol 
\pages 
\paperinfo  preprint (alg-geom/9609010)
\endref

\ref 
\key Gi1 
\by  A.B. Givental
\paper Homological geometry I.  Projective hypersurfaces
\jour Selecta Math.
\yr 1995
\vol 1
\pages 325--345
\endref

\ref 
\key  Gi2
\by A.B. Givental
\paper Homological geometry and mirror symmetry
\inbook Proc. Int. Congress of Math. I, Z\"urich 1994
\ed S.D. Chatterji
\yr 1995
\publ Birkh\"auser
\pages  472--480
\endref

\ref
\key  Gi3
\by A.B. Givental
\paper Equivariant Gromov-Witten invariants
\jour Internat. Math. Res. Notices
\yr 1996
\vol 13
\pages 1--63
\endref

\ref
\key  Gi4
\by A. Givental
\paper Stationary phase integrals, quantum Toda lattices, flag
manifolds and the mirror conjecture
\inbook Topics in Singularity Theory 
\bookinfo AMS Translations 180
\eds A. Khovanskii, A. Varchenko, and V. Vassiliev
\yr 1997
\publ Amer. Math. Soc.
\pages  103--115
\endref

\ref
\key  Gi5
\by A. Givental
\paper A mirror theorem for toric complete intersections
\inbook Topological field theory, 
primitive forms and related topics. 
\bookinfo  Proceedings of the 38th Taniguchi
Symposium, Progress in Math. 160
\eds  M. Kashiwara, A. Matsuo, K. Saito and I.
Satake
\yr 1998
\publ Birkh\"auser
\pages  141--175
\endref

\ref
\key Gi6
\by A. Givental
\paper A tutorial on quantum cohomology
\inbook Symplectic geometry and topology
\bookinfo IAS/Park City Math. Ser. 7
\publ Amer. Math. Soc.
\yr 1999
\pages 231--264
\endref

\ref 
\key  Gi-Ki
\by  A. Givental and B. Kim
\paper Quantum cohomology of flag manifolds and Toda lattices
\jour Commun. Math. Phys.
\yr 1995
\vol 168
\pages 609--641
\endref

\ref\key Gu1
\by M.A. Guest
\paper The topology of the space of rational curves on a 
toric variety
\jour Acta Math.
\vol 174
\yr 1995
\pages 119--145
\endref

\ref
\key Gu2
\by M.A. Guest
\book Harmonic Maps, Loop Groups, and Integrable Systems
\publ Cambridge Univ. Press
\yr 1997
\endref

\ref
\key  Gu-Ot
\by  M.A. Guest and T. Otofuji
\paper Quantum cohomology and the periodic Toda lattice
\paperinfo to appear 
\jour Commun. Math. Phys. (math.QA/9812127)
\yr 
\vol 
\pages 
\endref

\ref 
\key  Hi
\by F. Hirzebruch
\paper \"Uber eine Klasse von einfach-zusammenh\"angenden
komplexen Mannigfaltigkeiten
\jour Math. Ann.
\vol 124
\yr 1951
\pages 77--86
\endref

\ref
\key  Ki
\by B. Kim
\paper Quantum cohomology of flag manifolds $G/B$ and
quantum Toda lattices
\jour Ann. of Math.
\yr 1999
\vol 149
\pages 129--148
\endref

\ref
\key Ko
\by M. Kontsevich
\paper Enumeration of rational curves via torus actions
\inbook The Moduli Space of Curves
\bookinfo Progress in Mathematics 129
\publ Birkh\"auser
\yr 1995
\eds R. Dijkgraaf, C. Faber and G. van der Geer
\pages  335--368
\endref

\ref
\key  LVW
\by W. Lerche, C. Vafa and N.P. Warner
\paper Chiral rings in $N=2$ superconformal theories
\jour Nuclear Physics 
\yr 1989
\vol B324
\pages 427--474
\endref

\ref
\key Ma
\by Y.I. Manin
\book Frobenius Manifolds, Quantum Cohomology, and Moduli Spaces
\bookinfo Amer. Math. Soc. Colloquium Publications 47 
\publ Amer. Math. Soc.
\yr 1999
\endref

\ref
\key Mc-Sa
\by D. McDuff and D. Salamon
\book J-holomorphic Curves and Quantum Cohomology
\bookinfo Univ. Lec. Series 6
\publ Amer. Math. Soc.
\yr 1994
\endref

\ref
\key Od
\by T. Oda
\book Convex Bodies and Algebraic Geometry: An Introduction to the 
Theory of Toric Varieties
\publ Springer
\yr 1988
\endref

\ref
\key  RRW
\by M.S. Ravi, J. Rosenthal and X. Wang
\paper Degree of the generalized Pluecker embedding of a Quot scheme 
and quantum cohomology
\jour  Math. Ann.
\yr 1998
\vol 311
\pages 11--26
\endref

\ref
\key  Ru-Ti
\by Y. Ruan and G. Tian
\paper A mathematical theory of quantum cohomology
\jour J. Differential Geom.
\yr 1995
\vol 42
\pages 259--367
\endref

\ref
\key Si
\by B. Siebert
\paper An update on (small) quantum cohomology
\inbook Mirror symmetry, III (Montreal, PQ, 1995)
\bookinfo AMS/IP Stud. Adv. Math., 10
\publ Amer. Math. Soc.
\yr 1999
\pages 279--312
\endref

\ref
\key  Si-Ti
\by B. Siebert and G. Tian
\paper On quantum cohomology rings of Fano manifolds and
a formula of Vafa and Intriligator
\jour Asian J. Math.
\yr 1997 
\vol 1
\pages  679--695 (alg-geom/9403010)
\endref

\ref
\key Ti
\by G. Tian
\paper Quantum cohomology and its associativity
\inbook Current Developments in Mathematics, 1995
\publ Internat. Press
\yr 1994
\pages 361--401
\eds R. Bott {\it et al}
\endref

\ref
\key Vo
\by C. Voisin
\book Sym\acuteaccent etrie miroir
\bookinfo Panoramas et Synth\graveaccent eses 2
\publ Soc. Math. de France
\yr 1996
(English translation: 
Mirror Symmetry, SMF/AMS Texts and Monographs 1, 
Amer. Math. Soc. 1999)
\endref

\ref
\key We
\by A. Weinstein
\book Lectures on Symplectic Manifolds
\publ Amer. Math. Soc.
\bookinfo Regional Conference Series in Math. 29
\yr 1977
\endref

\endRefs
\enddocument